\documentclass[12pt]{icm}

\usepackage{amssymb, latexsym, graphicx}
\usepackage{color}
\usepackage{epic, eepic}

\newcommand{\E}{\mathcal E}
\newcommand{\PP}{\mathcal P}
\newcommand{\C}{\mathbb C}

\renewcommand{\phi}{\varphi}
\renewcommand{\epsilon}{\varepsilon}
\newcommand{\Z}{\mathcal Z}
\newcommand{\R}{\mathbb R}
\newcommand{\la}{\lambda}

\newtheorem{theorem}{Theorem}
\newtheorem{lemma}{Lemma}

\newtheorem{corollary}{Corollary}

\newtheorem{definition}{Definition}

\newtheorem{question}{Question}

\numberwithin{theorem}{section}
\numberwithin{claim}{section}
\numberwithin{question}{section}
\numberwithin{equation}{section}
\numberwithin{lemma}{section}
\numberwithin{corollary}{section}
\numberwithin{definition}{section}

\newcommand{\Di}{\operatorname{Di}}

\def\done{{1\hskip-2.5pt{\rm l}}}

\newcommand{\e}{\varepsilon}
\newcommand{\HH}{\mathcal H}
\newcommand{\Ss}{\mathbb S^2}

\def\done{{1\hskip-2.5pt{\rm l}}}

\title{Random Complex Zeroes and Random Nodal Lines}
\author[Fedor Nazarov and Mikhail Sodin]{Fedor Nazarov and Mikhail Sodin
\thanks{This work is partially supported by grant No. 2006136
of the United States - Israel Binational Science Foundation.
}
}

\contact[nazarov@math.wisc.edu]{F.N.:
Mathematics Department,
University of Wisconsin-Madison,
480 Lincoln Dr., Madison WI 53706
USA}
\contact[sodin@post.tau.ac.il]{M.S.: School of Mathematics, Tel Aviv University,
Tel Aviv 69978, Israel }

\begin{document}

\hfill{\sc To the memory of Oded Schramm}

\begin{abstract}

In these notes, we describe the recent progress in understanding the zero
sets of two remarkable Gaussian random functions: the Gaussian entire function
with invariant distribution of zeroes with respect to isometries of the complex
plane, and Gaussian spherical harmonics on the two-dimensional sphere.

\end{abstract}

\maketitle

\noindent These notes consist of two almost independent parts.
In both of them, we talk about zeroes
of special Gaussian random functions.
To understand them, we had to combine various
tools from complex and real analysis with
rudimentary probabilistic methods.
We think that the results and techniques presented here
can serve as guidelines
in other problems of similar nature arising in
analysis, mathematical physics, and probability theory.

The function $F$ that we consider in the first part is
a random analytic
function of one complex variable. In this case, one can recover
the zeroes of $F$ by applying the Laplacian to $\log |F|$.
This paves the way for using complex analysis tools, and
for this reason, the problems that we discuss in the first
part are pretty well understood by now,
though some intriguing questions still
remain open.

In the second part, we deal with topological properties
of the zero sets of random (real-valued) functions of several real variables.
This is an area with wealth of interesting and difficult questions and with
very few advances. In essence, in this part,
the reader will find a discussion of one recent theorem on
the number of connected components of zero sets of Gaussian
spherical harmonics along with various open questions.

\bigskip

\bigskip

\centerline{\sc\Large Part~I. Random complex zeroes}

\bigskip\par\noindent
The study of zeroes of random polynomials and random analytic functions has a long history.
It started with the pioneering works of Kac, Littlewood,
Offord, Rice, and Wiener, and was later continued
by Hammersley, Kahane, Maslova, and many others. The subject was revived
in the 1990's by several groups of researchers (Bogomolny-Bohigas-Leboeuf,
Shub-Smale, Edelman-Kostlan, Ibragimov-Zeitouni, Hannay, Bleher-Shiffman-Zelditch,
Nonnenmacher-Vo\-ros) who came from very different areas
and established new links to mathematical physics, probability
theory, and complex geometry. Some of these results were surveyed in
the lectures by Zelditch~\cite{Zelditch} and Sodin~\cite{Sodin-EMC};
see also an introductory article~\cite{NS0} and the
recent book by Hough, Krishnapur, Peres, and Vir\'ag~\cite{HKPV}.

In particular, Kostlan, Bogomolny-Bohigas-Leboeuf, Shub-Smale, and Hannay
introduced a remarkable construction of random Gaussian entire functions
with translation invariant distribution of their zeroes. Let
\[
F(z) = \sum_{n\ge 0} \zeta_n \frac{z^n}{\sqrt{n!}}
\]
where $\zeta_n$ are independent standard complex Gaussian random coefficients (i.e., the density
of $\zeta_k$ with respect to the Lebesgue measure in $\C$ is $ \tfrac1{\pi}e^{-|\zeta|^2}$).
The distribution of the random function $F$ is invariant with respect to rotations around
the origin, but it is {\em not} translation invariant,
for instance, because $\E |F(z)|^2 = e^{|z|^2}$ (here and below,
$\E$ means the mathematical expectation). However, the distribution of the zero set
$\Z=F^{-1}\{0\}$ {\em is} translation invariant. One of the ways to see this is to check
that the Gaussian random function
\[
F_\la (z) = F(z+\la) e^{-z\overline{\la} - \frac12|\la|^2}\,, \qquad \la\in\C\,,
\]
has the same covariance function as $F$:
\[
\E \bigl\{ F(z) \overline{F(w)} \bigr\}
= \E \bigl\{ F_\la(z) \overline{F_\la(w)} \bigr\} = e^{z\overline{w}}
\]
which is nothing else but the reproducing kernel in the classical Fock-Bargmann
space of entire functions. This coincidence is not accidental~\cite{NS0}.
Moreover, due to remarkable Calabi's rigidity~\cite[Section~2.5]{HKPV}, this is the
{\em only} translation invariant zero set of a Gaussian entire function
up to scaling.
We call the function $F$ the Gaussian Entire Function (G.E.F., for short).

It is worth mentioning that there exist similar constructions for other domains
with transitive groups of isometries (the hyperbolic plane, the Riemann sphere,
the cylinder and the torus).

\section{Linear statistics}\label{subsect:variance}

One of the most traditional ways to study asymptotic
properties of a random point process  $Z$ in the plane
is to take a test-function $h$, and to look at the asymptotic behaviour of
the linear statistics
\[
n_Z(r, h) = \sum_{a\in Z} h\bigl( \tfrac{a}{r} \bigr)
\]
as $r\to\infty$. We put $n(r, h)=n_\Z (r, h)$.
An easy computation shows that
\[
\E n(r, h) = \frac{r^2}{\pi} \int_{\R^2} h\,.
\]
If $h=\done_E$ is the indicator function of a set $E$, then
$n(r, \done_E) = n(rE)$ is the number of zeroes in the set $rE$.

\medskip
A usual ``triad'' in the study of linear statistics is

\medskip\centerline{\sc\small
variance, \quad asymptotic normality, \quad large fluctuations
}
\medskip\noindent
First, we'll discuss the variance, which is the easiest
part of the triad.

\subsection{The variance}

\begin{theorem}[The variance]\label{thm_variance}
For every non-zero function $h\in (L^1 \cap L^2)(\R^2)$ and every $R>0$,
\[
\operatorname{Var} n (r, h) = r^2\int_{\R^2}|\widehat h(\lambda)|^2 M(r^{-1}\lambda)\,{\rm d}m(\lambda)\,
\]
where
\[
M(\lambda)=\pi^3|\lambda|^4\sum_{\alpha\ge 1}\frac 1{\alpha^3}e^{-\frac{\pi^2}{\alpha}|\la|^2}\,,
\]
and
\[
\widehat{h}(\la) = \int_{\R^2} h(x) e^{-2\pi {\rm i}\, \langle \la, x \rangle }\, {\rm d}m(x)
\]
is the Fourier transform of $h$.
\end{theorem}
This theorem was proven in~\cite{NS1}.
The asymptotic of the variance had been known for two special cases since the work by Forrester and Honner~\cite{FH}:
if $h\in C^2_0$ (i.e., $h$ is a
$C^2$-function with compact support), then
\begin{equation}\label{eq:var1}
\operatorname{Var} n (r, h) = \frac{\zeta (3) + o(1)}{16\pi r^2}\, \|\Delta h\|_{L^2}^2\,,
\qquad  r\to\infty\,,
\end{equation}
while for bounded domains $G$ with piecewise smooth boundary,
\begin{equation}\label{eq:var2}
\operatorname{Var} n(rG) = \frac{\zeta(3/2)+o(1)}{8\pi^{3/2}}\, r\, L(\partial G)\,,
\qquad r\to\infty\,.
\end{equation}
Here, $\zeta(\,\cdot\,)$ is Riemann's zeta-function.

\medskip
A less precise form of Theorem~\ref{thm_variance} might be more illustrative:
\begin{equation}\label{eq_variance}
\operatorname{Var} n(r, h) \simeq r^{-2} \int_{|\lambda|\le r} |\widehat{h}(\lambda)|^2
|\la|^4 \, {\rm d}m(\la) + r^2
\int_{|\la|\ge r} |\widehat{h}(\la)|^2 \, {\rm d}m(\la) \,,
\end{equation}
where the notation $A\simeq B$ means that the quotient $B/A$ is bounded from below and from above by
positive numerical constants.
The right-hand side of~\eqref{eq_variance} interpolates $\| h \|^2_{L^2(m)}$ and
$\| \Delta h \|^2_{L^2(m)}$.

\subsection{Digression: ``superhomogeneous'' point processes}\mbox{}

\medskip\par\noindent
By~\eqref{eq:var2}, the random zero process $\Z$ belongs to the family of so called superhomogeneous
translation invariant point processes with fluctuations
of the number of points in large domains proportional to the length of the boundary rather than to the area, as it would be, say, for the Poisson process.

A ``toy model'' for such processes is the point process
\begin{equation}\label{eq:Gauss-perturb-lattice}
\mathcal S = \left\{ \omega + \zeta_\omega\colon \omega\in\sqrt{\pi}\,\mathbb Z^2 \right\}
\end{equation}
obtained by perturbing the lattice $\sqrt{\pi}\,\mathbb Z^2$
by independent standard complex Gaussian random variables
$\zeta_\omega$. The normalization by $\sqrt{\pi}$ is not essential here, it is introduced to
have asymptotically the same mean number of points in large areas as our process $\Z$ has.
The choice of the square lattice is not essential either.

Curiously, the same kernel $e^{z\overline w}$ that occurs in the definition of random
complex zeroes generates by a very different construction another interesting superhomogeneous
point process $\mathcal G$, namely, the determinantal process
whose $k$-point functions can be expresed in terms of the determinants formed by this kernel.
\[
\rho(z_1, ..., z_k) = \pi^{-k}e^{-\sum_{i=1}^k|z_i|^2}
\det \bigl\| e^{z_i \overline{z_j}} \bigr\|_{1\le i,j \le k}\,.
\]
This process arises as the large $N$ limit of eigenvalues of Ginibre ensemble of
$N\times N$ matrices with independent
standard complex Gaussian entries, and we will call it the {\em limiting Ginibre process}.
\begin{figure}[h]\label{fig:threeprocesses}
\centering
\includegraphics[height=1.4in]{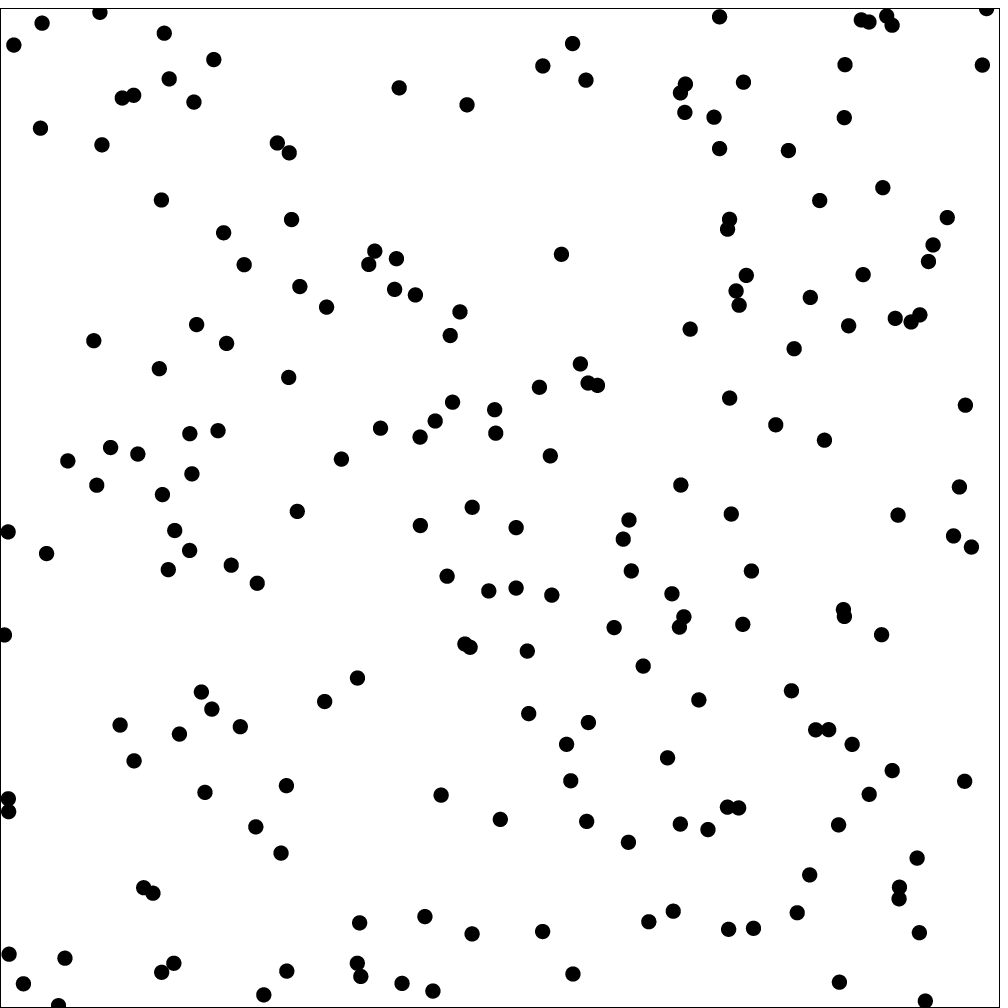}\hspace{.25in}
\includegraphics[height=1.4in]{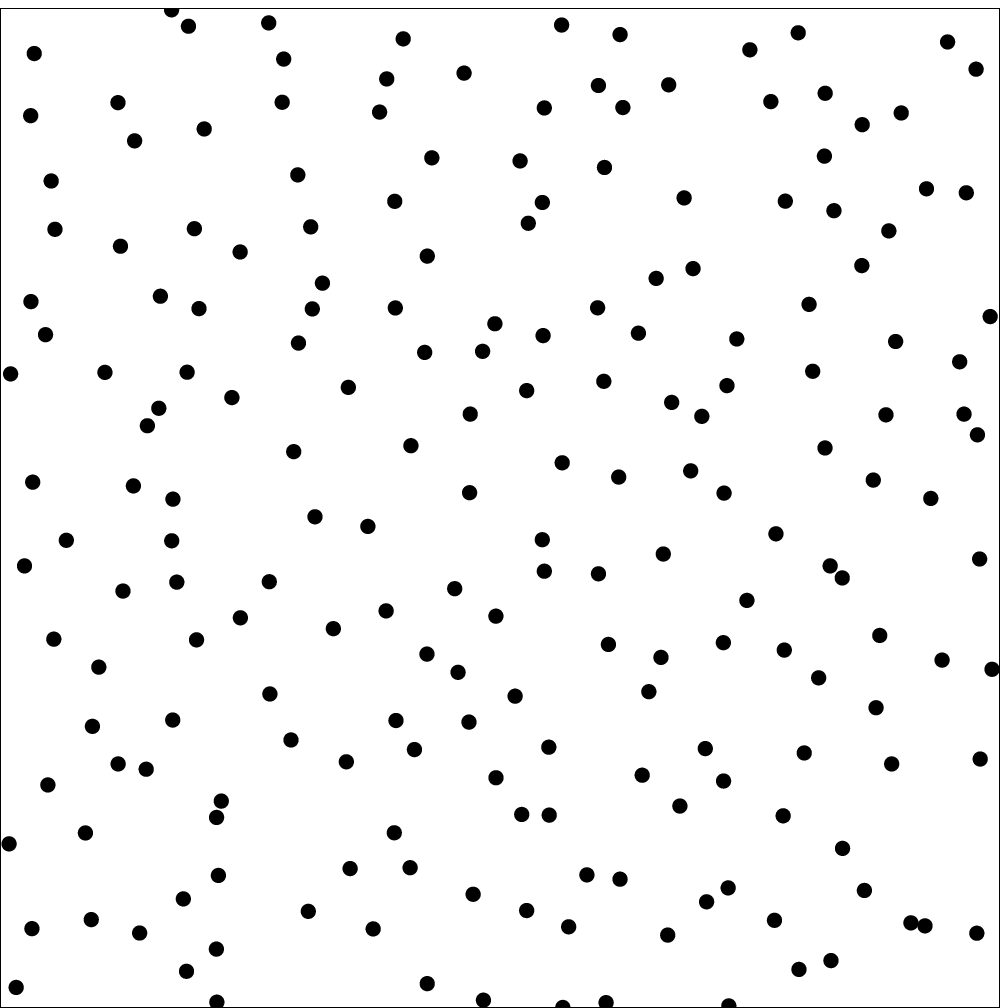}\hspace{.25in}
\includegraphics[height=1.4in]{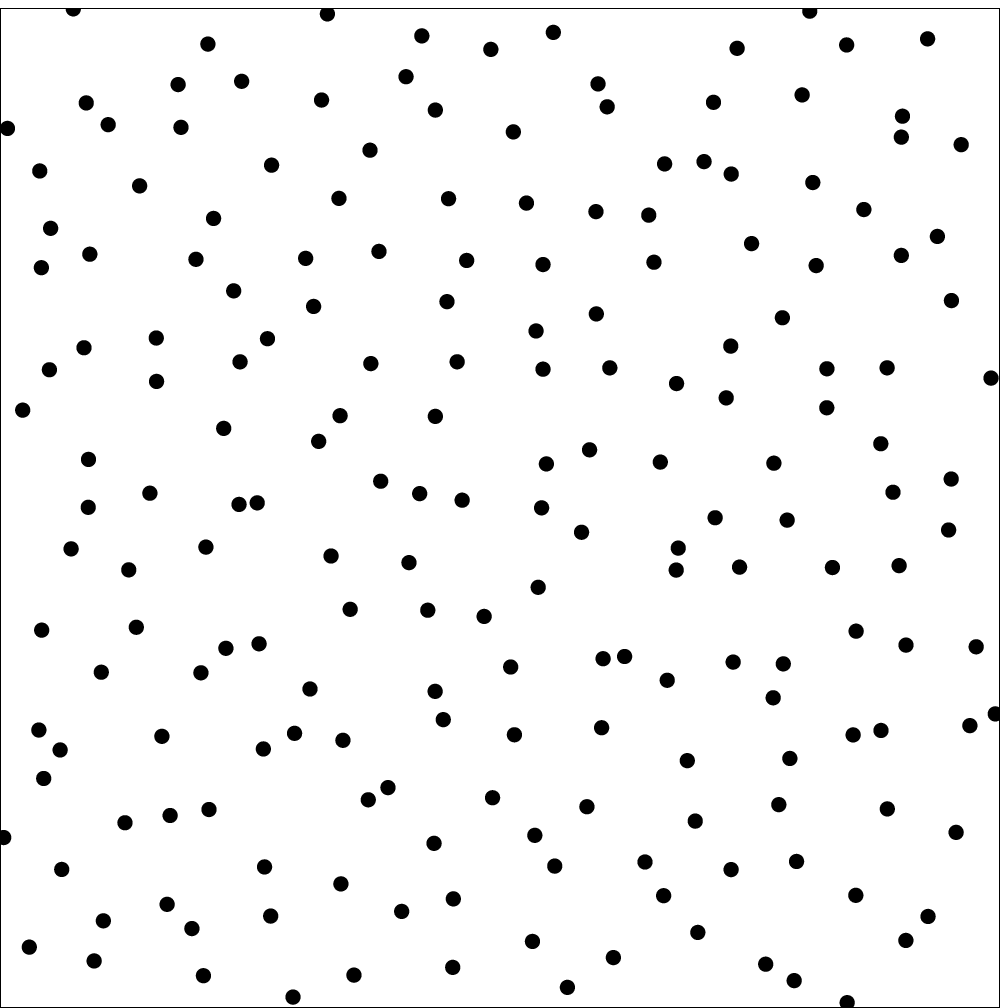}
\caption{\label{fig:gin&gaf}
Samples of the Poisson process (figure by B.~Vir\'{a}g),
limiting Ginibre process, and zeroes of a GEF (figures by M.~Krishnapur).
Some properties of the last two processes are quite different, though the eye does not easily
distinguish between them.
}
\end{figure}
It is known that the Ginibre point process is a special, explicitly solvable case of a one-component plasma
of charged particles of one sign confined by a uniform background of the opposite sign.
Though the one-component plasma has been studied by physicists for a long time,
it seems that most of rigorous results still pertain only to the very special case of
the Ginibre ensemble.

Resemblances and differences between the processes $\Z$, $\mathcal S$, and $\mathcal G$
were discussed both in the physical and the mathematical literature.
For instance, the behaviour of smooth linear statistics
for these three processes is quite different. In particular,
decay of the variance of smooth linear statistics~\eqref{eq:var1} distinguishes the zero
process $\Z$ from the processes $\mathcal G$ and $\mathcal S$, since for the latter two processes,
the variance of smooth linear statistics tends to the positive limit proportional to
$\| \nabla h \|^2_{L^2(m)}$.

\subsection{Asymptotic normality of fluctuations}\label{subsect:as-norm}

\subsubsection{Normal fluctuations}
We say that the linear statistics $n(r, h)$ have
{\em asymptotically normal fluctuations} if
the normalized linear statistics
\[
\frac{n(r, h) - \E n(r, h)}{\sqrt{\operatorname{Var} n(r, h)}}
\]
converge in distribution to the standard (real) Gaussian random variable as $r\!\to\!\infty$.

Let $C^\alpha_0$, $\alpha>0$, be the class of compactly supported $C^\alpha$-functions, by $C^0_0$
we denote the class of bounded compactly supported measurable functions.
\begin{theorem}[Asymptotic normality]\label{thm_normality}
Suppose that $h\in C^\alpha_0$ with some $\alpha \ge 0$, and that for some $\epsilon>0$
and for every sufficiently big $r$, we have
\begin{equation}\label{eq_ii}
\operatorname{Var} n(r, h)> r^{-2\alpha+\epsilon}\,.
\end{equation}
Then the linear statistics $n(r, h)$ have asymptotically normal fluctuations.
\end{theorem}
Note that by \eqref{eq_variance}, we always have $\operatorname{Var} n(r, h) \ge c(h) r^{-2}$
with positive $ c(h) $ independent of $r$. Hence, for $\alpha>1$, condition \eqref{eq_ii}
holds automatically, and we obtain the following
\begin{corollary}\label{cor_1}
Suppose that $h\in C^\alpha_0$ with $\alpha>1$.
Then the linear statistics
$n(r, h)$ have asymptotically normal fluctuations.
\end{corollary}
Using estimate~\eqref{eq_variance}, one can show
that for {\em any bounded measurable set $E$ of positive area},
$ \operatorname{Var} n(rE) \gtrsim r $, cf.~\eqref{eq:var2}. Hence,
\begin{corollary}\label{cor_2}
Let $E$ be a bounded measurable set of positive area. Then the number
of random complex zeroes $ n(rE)$ on the set $rE$ has asymptotically normal
fluctuations.
\end{corollary}

\subsubsection{\bf Abnormal fluctuations of linear statistics}
Do there exist $C_0^\alpha$-func\-tions $h$ with
abnormal asymptotic behaviour of linear statistics
$n(r, h)$? The answer is ``yes'', and the simplest example is provided by
the function $h_\alpha =|x|^\alpha\psi(x)$, where $\psi$ is a smooth cut-off
that equals $1$ in a neighborhood of the origin.
Clearly, $h_\alpha \in C^\alpha_0$ and it is not difficult to show that
$ \operatorname{Var} n(r, h_\alpha)  \simeq r^{-2\alpha}$. This
shows that Theorem~\ref{thm_normality} is sharp on a rough power scale.
The reason for the loss of asymptotic normality is that only a small neighbourhood
of the origin where $h_\alpha$ loses its smoothness
contributes to the variance of $n(r, h_\alpha)$.
This neighbourhood contains a bounded number of zeroes of $F$,
which is not consistent with the idea of normal fluctuations.

\subsubsection{Comments and questions}
Theorem~\ref{thm_normality} was preceded by a result of Sodin and Tsirelson~\cite[Part~I]{ST}.
Using the moment method and the diagrams, they showed that the fluctuations are
asymptotically normal provided that $h\in C_0^2$. Their technique works in several other
cases, for instance, when $h=\done_G$ is the indicator function of a bounded
domain $G$ with a piecewise smooth boundary. However, it seems very difficult
to adapt it for proving Theorem~\ref{thm_normality} in its full generality.

In the case $\alpha>0$, the proof of Theorem~\ref{thm_normality} uses
a classical idea of S.Bernstein to
approximate the random variable $n(R, h)$ by a sum of a large number of independent random variables with
negligible error. Such approximation becomes possible only after we separate the high and the
low frequencies in $h$.
In this approach, independence appears as a
consequence of the almost independence of the values of the G.E.F. at large distances, which we'll discuss
below in Section~\ref{subsect:Alm-Ind}.
We do not know whether asymptotic normality holds for all functions
$h\in C^1_0$, or whether the condition
$r^{2\alpha} \operatorname{Var} n(r, h)\to\infty$ is already sufficient for asymptotic
normality of linear statistics associated with a $C^\alpha_0$-function.
Also, we believe that the assertion of Theorem~\ref{thm_normality} can
be extended to functions  $h\in C^{\alpha}\cap L^2_0$ with $-1<\alpha<0$
but our current techniques seem insufficient to handle this case properly.

In the case $\alpha=0$, the proof uses a different idea which comes from
statistical mechanics. First, we show that $k$-point functions of the zero process
$\Z$ are clustering, see Section~\ref{subsect:unif-clust} for the precise
statement. Then, using clustering,
we estimate the cumulants of the random variable $n(r, h)$.

It is interesting to juxtapose Theorem~\ref{thm_normality} with what is known for the limiting
Ginibre process $\mathcal G$ described above.
For bounded compactly supported functions $h$,
a counterpart of Theorem~\ref{thm_normality} is a theorem of Soshnikov. In~\cite{Soshn},
he proved among other things that for arbitrary determinantal point processes, the fluctuations of
linear statistics associated with a compactly supported bounded positive function are normal
if the variance grows at least as a positive power of expectation as the intensity tends to infinity.
A counterpart of the limiting case $\alpha = 2$ in Theorem~\ref{thm_normality} (that is, of the result
from~\cite[Part~I]{ST}) was recently found by Rider and Vir\'ag in \cite{RV}.
They proved that the fluctuations for linear statistics of process $\mathcal G$
are normal
when the test function $h$ belongs to the Sobolev space $W_1^2$. It is not clear whether there is
any meaningful statement interpolating between
the theorems of Soshnikov and Rider and Vir\'ag. It can happen that
our Theorem~\ref{thm_normality} simply has no interesting counterpart for the process $\mathcal G$.
It is also worth mentioning that the proofs in the determinantal case are quite different from ours.
They are based on peculiar combinatorial identities for the cumulants of linear statistics that are a
special feature of determinantal point processes.

\subsection{Probability of large fluctuations}
\label{sect:LD}\mbox{}

\medskip\par\noindent
Now, we turn to the probability of exponentially rare events that, for some $r\gg 1$,
$ | n(r, h) - \E n(r, h) | $ is much bigger than $\sqrt{\operatorname{Var} (n(r, h)}$.
Mostly, we consider the case when $h$ is the indicator function of the unit disk
$\mathbb D$; i.e., we deal with the number $ n(r) $ of random zero points in the disk of
large radius $r$ centered at the origin. Recall that $\E n(r) = r^2$ and
$\E \bigl\{ (n(r)-r^2)^2 \bigr\} \sim cr$ for
$r\to\infty$ (with some $c>0$).
Hence, given $\alpha\ge \tfrac12$, we need to find the order of decay of the probability
$\PP \bigl\{ |n(r)-r^2|>r^\alpha  \bigr\}$.
\subsubsection{\bf Na\"{i}ve heuristics}
The aforementioned similarity between
the zero process $\Z$ and independent complex Gaussian perturbations $\mathcal S$
of the lattice $\sqrt{\pi}\,\mathbb Z^2$ helps to guess the correct
answer.

We fix the parameter $\nu > 0$, and consider the random point set
$ \mathcal S_\nu = \{ \omega +
\zeta_{\omega}\}_{\omega\in\mathbb Z^2} $, where $\zeta_\omega$
are independent, identical, radially distributed random variables
with the tails  $\PP \{ |\zeta_\omega|> t \} $ decaying as $\exp (-t^\nu)$
as $t\to\infty$. Set \[ n_{\nu}(r) = \# \{\omega\in \sqrt{\pi}\,\mathbb
Z^2\colon |\omega+\zeta_\omega|\le r \}.\]
Then, for every $ \alpha \ge\tfrac12 $ and every $\epsilon>0$,
\[
\exp[-r^{\phi(\alpha, \nu)+\epsilon}] < \PP \bigl\{ |n_\nu(r)-r^2|>r^\alpha  \bigr\}
< \exp[ -r^{\phi(\alpha, \nu)-\epsilon}]\,,
\]
provided that $r$ is sufficiently big. Here
\[
\phi (\alpha, \nu) =
\begin{cases}
2\alpha-1, & \tfrac12 \le \alpha \le 1; \\
(\nu+1)\alpha-\nu, & 1 \le \alpha \le 2; \\
(\tfrac12\nu +1)\alpha, & \alpha \ge 2\,.
\end{cases}
\]
Actually, one can find much sharper estimates for
$ \PP \bigl\{ |n_\nu(r)-r^2|>r^\alpha  \bigr\} $.

This suggests that the probability $ \PP \{ |n(r)-r^2|>r^\alpha \} $ we are after
should decay as $\exp [ - r^{\phi (\alpha) }]$ with
\[
\phi(\alpha) = \phi (\alpha, 2) =
\begin{cases}
2\alpha-1, & \tfrac12 \le \alpha \le 1; \\
3\alpha-2, & 1 \le \alpha \le 2; \\
2\alpha, & \alpha \ge 2\,.
\end{cases}
\]

\subsubsection{\bf Jancovici-Lebowitz-Manificat Law}
Unfortunately, we do not know how to represent random complex zeroes as independent,
or weakly correlated, Gaussian perturbations of the lattice points, so we cannot
use the heuristics given above. Nevertheless, we can prove
\begin{theorem}[JLM Law for random complex zeroes]\label{thm-JLM}
For every $ \alpha \ge\tfrac12 $ and every $\epsilon>0$,
\[
\exp [ -r^{\phi (\alpha) + \epsilon} ] < \PP\left\{ | n(r)-r^2 |
> r^\alpha \right\} < \exp[  -r^{\phi (\alpha) - \epsilon} ]
\]
for all sufficiently large $r>r_0(\alpha, \epsilon)$ with the same
$\phi (\alpha)$ as above.
\end{theorem}
In~\cite{JLM}, Jancovici, Lebowitz and Manificat showed that this law holds for the one-component
plasma. Their derivation was not a rigorous one, except for the case of
the limiting Ginibre process $\mathcal G$. It would be desirable to have
a clear explanation why {\em the same} Jancovici-Lebowitz-Manificat
law holds for the random processes $\Z$, $\mathcal S$, and $\mathcal G$ in the range $\alpha>1$.

\subsubsection{Comments and questions}
The function $\phi$ from the exponent in the Jancovici-Lebowitz-Manificat Law
loses smoothness at three points.
Accordingly, there are three different r\'egimes
($\frac12<\alpha<1$, $1<\alpha<2$, and $\alpha > 2$).
The point $\alpha = \tfrac12$ corresponds to the
asymptotic normality of $n(r)$, and deviations in the range $\tfrac12 < \alpha < 1$
are called {\em moderate}. In this range, the deviation $|n(r)-r^2|$ is small compared
to the length of the circumference $\{|z|=r\}$.
In this case, the theorem was proven by Nazarov, Sodin,
and Volberg~\cite{NSV2}.
The point $\alpha=1$ corresponds to the classical large deviations principle.
In the range $1<\alpha<2$, the deviation is already big compared to the length of
the boundary circumference, but is still small compared to the area of the
disk $\{|z|\le r\}$. In this case, the lower bound for $ \PP\left\{ | n(r)-r^2 |
> r^\alpha \right\}$ is due to Krishnapur~\cite{Krishnapur}, while the upper
bound was proven in~\cite{NSV2}.

The case $\alpha=2$ contains an estimate for the
``hole probability'' $ \PP \left\{ n(r)=0 \right\} $. In this case, the theorem was proved
by Sodin and Tsirelson~\cite[Part~III]{ST}. A very sharp estimate of the hole probability
\[
\log \PP  \left\{ n(r)=0 \right\}  = - \frac{3e^2}4 r^4 + O\bigl( r^{\frac{18}{5}} \bigr),
\qquad r\to\infty\,,
\]
was recently obtained by Nishry~\cite{Nishry}; in~\cite{Nishry2} he extended this asymptotics
to a rather wide class of entire functions represented by Gaussian Taylor series.
There are two interesting questions pertaining to the hole probability.
We have no idea how to find the asymptotics of the expected number of random complex
zeroes in the disk $R\mathbb D$, $R\ge r$, conditioned on the hole $\{ n(r)=0 \}$.
We also do not know how to extend Nishry's result from the unit disk to other
bounded domains $G$. It seems plausible that for a large class of bounded domains
$G$,
\[
\log \PP\left\{ n(r G) = 0 \right\} = -(\kappa (G) + o(1) ) r^4\,,  \qquad r\to\infty,
\]
with $\kappa (G)>0$.
If this is true, how does $\kappa (G)$ depend on $G$?

The range $\alpha>2$ in the Jancovici-Lebowitz-Manificat Law
is the ``overcrowding'' r\'egime. In~\cite{Krishnapur}, Krishnapur
proved that for $\alpha>2$,
\[
\log \PP \left\{ n(r) > r^\alpha \right\}
= - \bigl( \tfrac12\alpha -1 + o(1) \bigr) r^{2\alpha}\log r\,, \qquad
r\to\infty\,.
\]

The bounds in Theorem~\ref{thm-JLM} are not too tight. As we've already mentioned,
in some cases, much better bounds are known.
It would be good to improve precision of Theorem~\ref{thm-JLM} in other cases. For instance,
to show that for $\alpha\le 2$ and for $\delta>0$ there exists the limit
\[
\lim_{r\to\infty} \frac{\log  \PP\left\{ | n(r)-r^2 |
> \delta r^\alpha \right\}}{r^{\phi (\alpha)}}
\]
and to find its value.

\subsubsection{\bf Moderate deviations for smooth linear statistics} Here is a recent result of
Tsirelson~\cite{Tsirel}:
\begin{theorem}\label{thm:tsirel}
Let $h\in C^2_0$. Then
\[
\log \mathbb P \bigl\{ rn(r, h) > t \sigma \|\Delta h\|_{L^2} \bigr\}
=(1+o(1)) \log \Bigl( \frac1{\sqrt{2\pi}} \int_t^\infty e^{-x^2/2}\, {\rm d}x \Bigr)
\]
and
\[
\log \mathbb P \bigl\{ rn(r, h) <- t \sigma \|\Delta h\|_{L^2} \bigr\}
=(1+o(1)) \log \Bigl( \frac1{\sqrt{2\pi}} \int_t^\infty e^{-x^2/2}\, {\rm d}x \Bigr)\,,
\]
as $r\to\infty$, $t>0$, and $t \tfrac{\log^2 r}{r}\to 0$. Here,
$\sigma^2 = \frac{\zeta (3)}{16\pi}$ (cf.~\eqref{eq:var1}).
\end{theorem}

The proof of this theorem is quite intricate.
Note that it gives bounds that are much sharper than the ones in Theorem~\ref{thm-JLM}.
In the case $t={\rm const}$, Theorem~\ref{thm:tsirel} gives another proof
of the asymptotic normality of smooth linear statistics of random complex zeroes.

It is not clear whether the assumption $t \tfrac{\log^2 r}{r}\to 0$ can be replaced by a more
natural one $ \tfrac{t}{r} \to 0 $.
To the best of our knowledge, until now, there have been no results about large
or huge deviations
for smooth linear statistics of random complex zeroes when $t$ is comparable or
much larger than $r$.

\section{Uniformity of spreading of random complex zeroes over the plane}

Let $Z$ be a point process in $\R^d$ with the
distribution invariant with respect to the isometries of $\R^d$.
A natural way to check how evenly the process $Z$ is spread over
$\R^d$ is to find out how far
the counting measure
\[ n_Z = \sum_{a\in Z} \delta_a \] of
the set $Z$ ($\delta_a$ is the unit mass at $a$) is from the Lebesgue
measure $m_d$ in $\R^d$. We describe a  convenient way to measure the distance between
$n_Z$ and $m_d$.

Suppose that the mean
number of points of $Z$ per unit volume equals $1$.
We want to partition the whole space $\R^d$, except possibly a subset of
zero Lebesgue measure, into disjoint sets $B(a)$ of Lebesgue measure
$1$ indexed by sites $a\in Z$
in such a way that each set $B(a)$ is located not too far from the corresponding site
$a\in Z$. In other words, we are looking for a measurable map $T\colon \R^d\to Z$
such that for each $a\in Z$, we have $m_d (T^{-1}\left\{ a \right\})=1$.
We also want the distances $ |Tx-x| $ to be not too large.
The map $T$ is called the {\em transportation}
(a.k.a.  ``matching'', ``allocation'',
``marriage'', etc.) of the Lebesgue measure $m_d$ to the set $Z$.

Alternatively, we can fix a lattice
$\Gamma\subset\R^d$ with cells of unit volume, and
look for a bijection $\Theta\colon \Gamma\to Z$ for which
the distances $|\Theta\gamma-\gamma|$, $\gamma\in\Gamma$, are not
too large. Since for each two lattices $\Gamma_1$, $\Gamma_2$ with cells
of the same volume, there is a bijection $\theta\colon \Gamma_1\to\Gamma_2$
with $ \sup \left\{ | \theta \gamma - \gamma |\colon \gamma\in\Gamma_1 \right\}<\infty $,
the choice of the lattice is not important, so we can take $\Gamma=\mathbb Z^d$.

Since we deal with {\em random} discrete sets $Z$,
the corresponding transportation maps $T$ (or the bijections $\Theta$)
will be {\em random maps}.
In interesting cases (including the random complex zeroes $\Z$),
almost surely, the transportation
distances $|Tx-x|$ are unbounded, so we are
interested in the rate of decay of the probability tails
$\PP \{ |Tx-x|>R \}$ as $R\to\infty$.

Here we present two approaches to this problem developed in~\cite[Part~II]{ST}
and in~\cite{NSV1}. Though we discuss only
the random complex zeroes $\Z$, we  believe that both approaches should work
for other natural translation invariant point processes.
At last, we recall that the random complex zero process $\Z$ has intensity $\pi$,
not $1$. For this reason, we will look for a transportation of the measure
$\pi m_2$ to $\Z$, and for a bijection between the lattice
$\sqrt{\pi}\,\mathbb Z^2$ and $\Z$.

\subsection{Random complex zeroes as randomly perturbed lattice points}\mbox{}

\begin{theorem}[Existence of well-localized bijection]\label{thm:existence}
There exists a translation invariant random function $\xi\colon \mathbb Z^2\to\C$
such that

\noindent\textup{(a)} the random set $ \{ \gamma  +
\xi (\gamma)\colon \gamma\!\in\!\sqrt{\pi}\,\mathbb Z^2 \} $ is equidistributed with the
random complex zeroes $\mathcal Z$;

\noindent\textup{(b)} $ \PP \{|\xi (0)| > R \} \le \exp \left(
-cR^4/\log R\right) $ for some $ c>0 $ and every $ R\ge 2 $.
\end{theorem}

The theorem is almost optimal since
the probability that the disk of radius $\la\ge 1$ is
free of random complex zeroes is not less than $\exp \left( -C\la^4\right)$.
It seems that the question about the
existence of a matching between  the lattice and $\mathcal Z$  with tails decaying as
$ \exp \left( - c\la^4\right) $ remains open as well as the same question for the
Gaussian perturbations $\mathcal S$ of the lattice points and for the limiting
Ginibre process $\mathcal G$.

It would be interesting to find a version of Theorem~\ref{thm:existence} with weakly correlated
perturbations $\xi_{k,l}$ at large distances. This could shed some light to the reasons hidden behind
the Jancovici-Lebowitz-Manificat Law.

\subsubsection{\bf Uniformly spread sequences in $\R^d$}
The proof of Theorem~\ref{thm:existence} is based on a deterministic idea
which might be useful in study of the uniformity of spreading of sequences and
measures.
We need to establish the bijection
between the sets $\mathcal Z$ and $\sqrt{\pi} \mathbb Z^2$ with controlled tails of the
distances $|\xi_{k, l}|$. First we look at a simpler situation when $|\xi_{k, l}|$ are uniformly
bounded. It is too much to expect from a typical zero set, but let us try anyway.
We say that the set $Z\in \R^d$ is {\em $r$-uniformly spread} over $\R^d$
(with density $1$) if there exists a bijection between $Z$ and a lattice with the unit volume of the
cell such that the distances between $Z$ and the corresponding lattice points do not exceed $r$.
If such a bijection exists then clearly
\begin{equation}\label{eq:unif-spread1}
n(U) \le \nu (U_{+r}) \qquad {\rm and} \qquad \nu (U)\le n(U_{+r})
\end{equation}
for every $U\subset\R^d$; here $U_{+r}$ stands for the $r$-neighbourhood of $U$, $n$ is the counting
measure of the set $Z$, and $\nu$ is the counting measure of
the lattice. In fact, \eqref{eq:unif-spread1} is not only necessary but also {\em sufficient}, which is
basically a well-known locally finite marriage lemma due to M.~Hall and R.~Rado.
When verifying condition~\eqref{eq:unif-spread1},
we can replace $\nu$ by the Lebesgue measure
$m_d$ at the expense of adding a constant to~$r$ .

Now, given a locally finite measure $\mu$ on $\R^d$, we define
$\Di (\mu)$ as the infimum of $r\in (0, \infty)$ such that
\[
\mu (X) \le m_d (X_{+r}) \qquad {\rm and} \qquad m_d (X)\le \mu(X_{+r})
\]
for every bounded Borel set $X\subset\R^d$. The range of $\Di$ is $[0, +\infty]$ with the
both ends included.
The following  theorem gives a useful upper bound for $\Di (\mu)$ in terms of the
potential $u$:
\begin{theorem}[Upper bound for the transportation distance]\label{thm:unif-spread}
Let $u$ be a locally integrable function in $\R^d$ such
that $\Delta u = \mu-m_d$ in the sense of distributions. Then
\[
\Di(\mu) \le \operatorname{Const}_d\cdot \inf_{r>0} \big\{
r + \sqrt{\|u*\chi_r\|_\infty} \, \big\} \, .
\]
Here, $ \chi_r $ is the indicator function of the ball of radius $r$
centered at the origin normalized by the condition $ \| \chi_r \|_{L^1} = 1$,
and $*$ denotes the convolution.
\end{theorem}

Now, we explain how Theorem~\ref{thm:existence} is deduced from
Theorem~\ref{thm:unif-spread}. After smoothing, the random potential
$U(z)= \log |F(z)| - \tfrac12 |z|^2$ is locally uniformly bounded.
Still, a.s. it remains unbounded in $\C$, so we cannot apply
Theorem~\ref{thm:unif-spread} directly. The idea is to introduce on $\C$ a
{\em random metric} $\rho$ that depends on a G.E.F. $F$. The metric $\rho$ is
small where the random potential $U$ is large. Then we apply a counterpart of
Theorem~\ref{thm:unif-spread} with the distances measured in the metric $\rho$,
instead of the Euclidean one.

\subsubsection{Comments}
Theorems~\ref{thm:existence} and~\ref{thm:unif-spread}
are taken from Sodin and Tsirelson~\cite[Part~II]{ST} (cf.~\cite{ST2}).
In that paper, the authors proved a weaker subgaussian estimate for the tails,
however, after a minor modification of the proof given therein, one gets
the result formulated here.
Note that the method developed in Sodin and Tsirelson~\cite[Part~II]{ST}
needs only the existence of a stationary random vector field $v$ in $\R^d$ with
$\operatorname{div}{v} = \mu-c_d m_d$. The tail estimate depends on the rate of decay of
the tails of the field $v$ or of the tails of the potential $u$ such that
$v=\nabla u$ (if such a $u$ exists).

In the last 20 years, the concept of {\em uniformly spread} discrete subsets of $\R^d$
has appeared in very different settings. Laczkovich used uniformly spread
sets in $\R^d$ in his celebrated
solution of the Tarski's circle squaring problem~\cite{Laczkovich1}
(see also~\cite{Laczkovich2}).
There are various probabilistic counterparts of this notion.
For instance, Ajtai, Koml\'os and Tusn\'ady~\cite{AKT}, Leighton
and Shor~\cite{LS}, and Talagrand~\cite{T}
studied a finite counterpart of this, namely, a high probability
matching of a system of $N^2$ independent random points in the square $[0, N)^2\subset
\R^2$ with the grid $\Z^2\cap [0, N)^2$.

\subsection{Gradient transportation}\label{section:GradTransp}\mbox{}

\medskip\par\noindent
Unfortunately, the proof of Theorem~\ref{thm:existence} is a pure existence one. It
gives us no idea about what the (almost) optimal transportation of the Lebesgue measure
to the zero process $\Z$ looks like. Now, we discuss another approach, namely, the
transportation by the gradient flow of a random potential.
The main advantage of this approach is that it yields a quite
natural and explicit construction for the desired transportation.

\subsubsection{\bf Basins of zeroes}
Let $U(z)=\log|F(z)|-\frac12{|z|^2}$ be the random potential
corresponding to the G.E.F. $F$. It is easy to check that the distribution of
$U$ is invariant with respect to the isometries of the plane.
We shall call any integral curve of the differential equation
\[
\frac{dZ}{dt} = -\nabla U(Z)
\]
a {\em gradient curve} of the potential $U$.
We orient the gradient curves in the direction of decrease of $U$
(this is the reason for our choice of the minus sign in the
differential equation above). If $z\notin \Z$, and
$\nabla U(z)\ne 0$, by $\Gamma_z$ we denote the (unique) gradient
curve that passes through the point $z$.

\begin{definition}{\rm Let $a$ be a zero of the G.E.F. $F$.
The {\it basin} of $a$ is the set
$$
B(a)=\{z\in \C\setminus \Z\colon \,\nabla U(z)\ne 0,\text{ and } \Gamma_z
\text{ terminates at }a\}\,.
$$
}
\end{definition}

Clearly, each basin $B(a)$ is a connected open set, and $B(a')\cap
B(a'')=\varnothing$ if $a'$ and $a''$ are two different zeroes of
$F$. Remarkably, {\em all bounded basins  have the same area $\pi$}.
Indeed, $\tfrac{\partial U}{\partial n} = 0$ on $\partial B(a)$
and therefore, applying the Green formula and recalling that the
distributional Laplacian of $U$ equals $\Delta U = 2\pi
\sum_{a\in\mathcal Z_F} \delta_a - 2m$, one gets
$$
1-\frac {m B(a)}{\pi}=\frac 1{2\pi}\iint_{B(a)}\Delta
U(z)\,{\rm d}m (z) = \frac1{2\pi}\int_{\partial B(a)}\frac {\partial
U}{\partial n}(z)\,|dz|=0\,;
$$
i.e., $m B(a) = \pi$.
The picture below helps to visualize what's going on.
\begin{figure}[h]\label{fig_manju}
\begin{center}
\scalebox{0.25}{\includegraphics{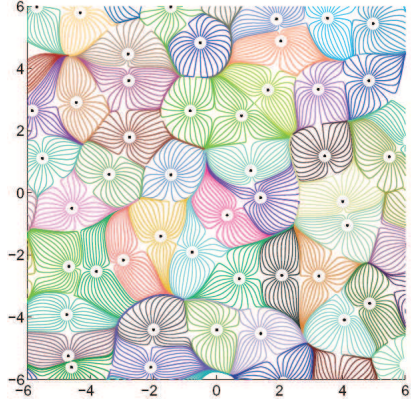}}
\end{center}
\caption{Random partition of the plane into domains of equal area generated
by the gradient flow of the random potential $ U $ (figure by M.~Krishnapur).
The lines are gradient
curves of $U$, the black dots are random zeroes. Many basins meet
at the same local maximum, so that two of them meet tangentially, while
the others approach it cuspidally forming long,
thin tentacles.}
\end{figure}

\subsubsection{\bf Results}
\begin{theorem}[Random partition]\label{thm:random-partition}
Almost surely, each basin is bounded by finitely many smooth
gradient curves (and, thereby, has area $\pi$), and
\[
\mathbb C=\bigcup_{a\in\Z} B(a)
\]
up to a set of measure $0$ (more precisely, up to countably many
smooth boundary curves).
\end{theorem}

The tails of this random partition have three characteristic exponents
${\large\bf 1}$, ${\large\bf \frac85}$, and ${\large\bf 4}$.
The probability that the diameter of a
particular basin is greater than $R$ is exponentially small in
$R$. Curiously enough, the probability that a given point $z$ lies at a distance
larger than $R$ from the zero of $F$ it is attracted to decays much faster:
as $e^{-R^{8/5}}$. This is related to long thin tentacles seen on the picture
around some basins. They increase the typical diameter of the basins though the probability
that a given point $z$ lies in such a tentacle is very small.
At last,  given $\epsilon>0$, the probability that it is impossible to throw away
$\epsilon\%$ of the area of the basin so that the diameter of the remaining part is less than
$R$ decays as $e^{-R^4}$. All three exponents are optimal.
The proofs of these results rely on the following
long gradient curve theorem.
\begin{theorem}[Long gradient curve]\label{thm:long-grad-curve}
Let $ R\ge 1 $. Let $Q(R)$  be the square centered at the origin with side length $R$.
The probability of the event that there exists a
gradient curve joining $\partial Q(R)$ with $\partial Q(2R)$ does not exceed
$Ce^{-cR(\log R)^{3/2}}$.
\end{theorem}
The proof of this theorem is, unfortunately, rather long and complicated.
It might be helpful for the reader to look at the first version of~\cite{NSV1}
posted in the {\tt arxiv}
where the authors gave a more transparent proof of a weaker upper bound
$Ce^{-cR\sqrt{\log R}}$ in the long gradient curve
theorem.

\subsubsection{Comments and questions}
Gradient transportation was introduced by Sodin and Tsirelson~\cite[Part~II]{ST}
and studied by Nazarov, Sodin, Volberg in~\cite{NSV1}.

There are several questions related to the statistics
of our random partition of the plane.
It is not difficult to show that, almost surely, any given
point $z\in\C$ belongs to some basin. We denote that basin by $B_z$,
and the corresponding sink by $a_z$.
We say that two basins are
neighbours if they have a common gradient curve on the boundary.
By $N_z$ we denote the number of basins $B$ neighbouring the basin
$B_z$. Clearly, $N_z$ equals the number of saddle points of the
potential $U$ connected with the sink $a_z$ by gradient curves.
Heuristically, since almost surely each saddle point is connected
with two sinks,
\[
\E N_z = 2\, \frac{{\rm mean\ number\ of\ saddle\ points\
per\ unit\ area}}{{\rm mean\ number\ of\ zeroes\ per\ unit\
area}}\,.
\]
Douglas, Shiffman and Zelditch proved in \cite{DSZ} that the mean
number of saddle points of $U$ per unit area is $\frac4{3\pi}$.
(They proved this for another closely related ``elliptic model''
of Gaussian polynomials. It seems that their proof also works for
G.E.F.'s) This suggests that $\E N_z = \frac83$.
Another characteristic of the random partition is the number
of basins that meet at the same local maximum. Taking into account
the result from \cite{DSZ}, we expect that its average equals $8$.

The next question concerns the topology of our random partition
of the plane. By the {\em skeleton} of the gradient flow we mean the connected
planar graph with vertices at local maxima of $U$ and edges
corresponding to the boundary curves of the basins. The graph may
have multiple edges and loops. We do not know
whether there are any non-trivial topological restrictions on finite parts
of the skeleton that hold almost surely.

\smallskip
In~\cite{CPPR} Chatterjee, Peled, Peres, Romik applied the ideas
from~\cite{NSV1} to study
the gradient transportation of the Lebesgue measure to the Poisson point
process in $\R^d$ with $d\ge 3$ (they called it `gravitational allocation').
Their work required a delicate and thorough analysis of the behaviour
of the Newtonian  potential of the Poisson point process.
It's worth mentioning that a very different construction of
the {\em stable marriage} between the Lebesgue
measure $m_d$ and the Poisson process in $\R^d$ with $d\ge 2$
was developed by Hoffman, Holroyd and Peres in~\cite{HHP}.
The case $d=2$ is especially interesting: see the recent work
by Holroyd, Pemantle, Peres, Schramm~\cite{HPPS}.

\section{Almost independence and correlations}

\subsection{Almost independence at large distances}\label{subsect:Alm-Ind}\mbox{}

\medskip\par\noindent
The covariance function
of the normalized Gaussian process $F^*(z)=F(z)e^{-\frac12 |z|^2}$ equals
\[ e^{z\overline{w}-\frac12|z|^2-\frac12 |w|^2} = e^{{\rm i}{\rm Im}(z\overline{w}) - \frac12 |z-w|^2}\,, \]
which decays very fast as $|z-w|$ grows. This suggests an idea that
the zeroes of G.E.F.'s must be ``almost independent'' on large distances.
Still the precise formulation of this independence property
is not obvious: due to analyticity
of $F$, if we know the process $F^*$ in a neighbourhood of some point, we know it everywhere on the plane.

It is not difficult to show that two standard complex Gaussian random variables
with small covariance can be represented as small perturbation of two
{\em independent} standard complex Gaussian random variables.
Developing this idea, we show that
if $\{K_j\}$ is a collection
of well-separated compact sets, then the restrictions $F^* \big|_{K_j}$ of normalized process $F^*$
can be simultaneously approximated by restrictions $F_j^* \big|_{K_j}$
of normalized {\em independent} realizations of G.E.F.'s $F_j$
with high precision and the probability very close to $1$.
This is a very useful principle that lies in the core of the proofs
of most of the results described above. Here is the precise statement~\cite{NS1}:
\begin{theorem}[Almost independence]\label{thm_almost_indep} Let $F$ be a G.E.F..
There exists a numerical constant $A>1$ with the
following property. Given a family of compact sets $K_j$ in $\C$ with diameters
$d(K_j)$, let $\la_j\ge \max \{d(K_j), 2 \}$.
Suppose that $A\sqrt{\log\la_j}$-neighbourhoods of the sets $K_j$ are pairwise disjoint. Then
\[
F^* = F_j^* + G_j^* \qquad  on \ K_j,
\]
where $F_j$ are independent G.E.F.'s and for every $j$, we have
\[
\PP \bigl\{ \max_{K_j} |G_j^*| \ge \la_j^{-1} \bigr\}
\lesssim  e^{-\la_j} \,.
\]
\end{theorem}
Less general versions of this result were proven in~\cite{NSV1, NSV2}.

The proof of Theorem~\ref{thm_almost_indep} goes as follows. First, for each compact set $K_j$,
we choose a sufficiently dense net $Z_j$ and consider the bunch $N_j = \bigl\{ v_z\colon z\in Z_j \bigr\}$
of unit vectors $v_z = F^*(z)$ in the Hilbert space of
complex Gaussian random variables. Since the compact sets
$K_j$ are well-separated, the bunches $N_j$ are almost orthogonal to each other. Then
we slightly perturb the vectors $v_z$ without changing the angles between
the vectors within each bunch $N_j$, making the bunches orthogonal to each other.
More accurately, we construct new bunches
$\widetilde{N}_j = \bigl\{ \widetilde{v}_z\colon z\in Z_j \bigr\}$ so that for $z\in Z_j$,
$\zeta\in Z_k$,
\[
\langle \widetilde{v}_z, \widetilde{v}_\zeta\rangle
=
\begin{cases}
\langle v_z, v_\zeta\rangle & {\rm for\ } j=k, \\
0 & {\rm for\ } j\ne k
\end{cases}
\]
with good control of the  errors
$ \| v_z - \widetilde{v}_z\| $.
Then we extend the Gaussian bunches
$\bigl\{ \widetilde{v}_z e^{\frac12 |z|^2}\colon z\in Z_j \bigr\}$
to {\em independent} G.E.F.'s $F_j$. The difference $G_j = F-F_j$ is a
random entire function
that is small on the net $Z_j$ with probability very close to one.
At the last step of the proof, using some simple complex analysis,
we show that $G_j^*$ is small everywhere on $K_j$.

\subsection{Uniform estimates of $k$-point functions. Clustering}\label{subsect:unif-clust}\mbox{}

\medskip\par\noindent
There is yet another way (originated in statistical mechanics)
to describe point processes by the properties of their $k$-point correlation functions. Recall that
the $k$-point function $\rho=\rho_k$ of the zero process $\Z$
is a symmetric function on $\C^k$ defined outside of the diagonal subset
\[ {\tt Diag}(\C^k) = \{(z_1, ..., z_k)\colon z_i= z_j {\rm \ for\ some\ } i\ne j \} \]
by the formula
\begin{equation}\label{eq_k-pt}
\rho (z_1, ..., z_k) =
\lim_{\epsilon\to 0} \frac{p_\epsilon (z_1, ..., z_k)}{(\pi\epsilon^2)^k}
\end{equation}
where $p_\epsilon (z_1, ..., z_k)$ is the probability that each disk
$\{|z-z_j|\le \epsilon \}$, $1\le j \le k$, contains at least one point of
$\Z$.
The $k$-point functions describe correlations within $k$-point
subsets of the point process. Estimates for the $k$-point functions are crucial
for understanding many properties of point processes.
The following results taken from~\cite{NS2}
provide rather complete quantitative information about the
behaviour of the $k$-point functions of random complex zeroes.

The first result treats the local behaviour of $k$-point functions.
It appears that for a wide class of non-degenerate Gaussian
analytic functions, the $k$-point functions of their zeroes
exhibit universal local repulsion when some of the variables
$z_1, ..., z_k$ approach each other.

Recall that a Gaussian analytic function $ f(z) $ in a plane domain
$ G\subseteq \C $ is the sum
\[
f(z) = \sum_n \zeta_n f_n (z)
\]
of analytic functions $ f_n(z) $ such that
\[
\sum_n |f_n(z)|^2 < \infty \quad \text{locally  uniformly  on } G,
\]
where $\zeta_n$ are independent standard complex Gaussian coefficients.
By $\rho_f = \rho_f(z_1, ..., z_k)$ we denote the $k$-point function of the zero
set of the function $f$. It is a symmetric function defined outside the diagonal
set $ {\tt Diag}(G^k) $ as in~\eqref{eq_k-pt}.

We skip the technical definition of $d$-degeneracy, which we use
in the assumptions of the next theorem, and only mention
that Gaussian Taylor series (either infinite, or finite)
\[
f(z) = \sum_{n\ge 0} \zeta_n c_n z^n
\]
are $d$-nondegenerate,  provided that $c_0, c_1, ..., c_{d-1} \ne 0$.
In particular, the G.E.F. is $d$-nondegenerate for every positive integer $d$.
\begin{theorem}[Local universality of repulsion]\label{thm-weak-univ}
Let $f$ be a $2k$-nondegenerate Gaussian analytic function
in a domain $G$, let $\rho_f $ be a $k$-point
function of zeroes of $f$, and let $K\subset G$ be a compact set.
Then there exists a positive constant $C=C(k, f, K)$ such that,
for any configuration of pairwise distinct points $z_1, ..., z_k\in K$,
\[
C^{-1} \prod_{i<j} |z_i-z_j|^2 \le \rho_f (z_1, ..., z_k)
\le C \prod_{i<j} |z_i-z_j|^2\,.
\]
\end{theorem}

The next result is a clustering property of zeroes of G.E.F.'s.
It says that if the
variables in $\C^k$ can be split into two groups located far from each other,
then the function $\rho_k$ almost equals the product of the corresponding
factors. This property is another manifestation of almost independence of
points of the process at large distances. It plays a central r\^{o}le in the
proof of the asymptotic normality theorem~\ref{thm_normality} for bounded
measurable functions.

For a non-empty subset $I = \left\{ i_1, ..., i_\ell \right\} \subset \left\{ 1, 2, ..., k \right\}$, we
set $Z_I = \left\{ z_{i_1}, ..., z_{i_\ell} \right\}$. We denote by
\[ d(Z_I, Z_J) = \inf_{i\in I, j\in J} |z_i-z_j| \] the distance
between the configurations $ Z_I $ and $ Z_J $.
\begin{theorem}[Clustering property]\label{thm_main}
For each $k\ge 2$, there exist positive constants $C_k$ and $\Delta_k$ such that for
each configuration $Z$ of size $k$ and each partition of the set of indices
$\left\{ 1, 2, ..., k \right\}$ into two non-empty
subsets $ I $ and $ J $ with $ d(Z_I, Z_J) \ge 2\Delta_k $, one has
\begin{equation}\label{eq:cluster}
1 - \epsilon
\le \frac{ \rho (Z) }{ \rho (Z_I) \rho (Z_J)}
\le 1 + \epsilon
\quad{\rm with} \quad
\epsilon =  C_k e^{-\frac12 (d(Z_I, Z_J)-\Delta_k)^2}\,.
\end{equation}
\end{theorem}

Combining Theorems~\ref{thm-weak-univ} and~\ref{thm_main}, and taking into
account the translation invariance of the point process $\Z$, we obtain
a uniform estimate for $\rho_k$ valid in the whole $\C^k$:
\begin{theorem}\label{thm-multipl} For each $k\ge 1$,
there exists a positive constant $C_k$ such that for each configuration
$(z_1, ..., z_k)$,
\[
C_k^{-1} \prod_{i<j} \ell (|z_i-z_j|)
\le \rho (z_1, ..., z_k) \le C_k \prod_{i<j} \ell(|z_i-z_j|)\,,
\]
where $\ell (t) = \min(t^2, 1)$.
\end{theorem}

The  proofs of Theorems~\ref{thm-weak-univ} and~\ref{thm_main}
start with the classical Kac-Rice-Hammersley  formula~\cite[Chapter~3]{HKPV}:
\begin{equation}\label{eq:KRH}
\rho_f (z_1, ..., z_k)
= \int_{\C^k} |\eta_1|^2 ... |\eta_k|^2\, \mathcal D_f(\eta'; z_1, ..., z_k)\,
{\rm d} m(\eta_1) ... {\rm d} m(\eta_k),
\end{equation}
where $\mathcal D_f(\, \cdot \,; z_1, ..., z_k)$
is the density of the joint probability distribution of the random variables
\begin{equation}\label{eq-20}
f(z_1),\, f'(z_1),\, ...\,, f(z_k),\, f'(z_k)\,,
\end{equation}
and $\eta' = \left( 0, \eta_1, ..., 0, \eta_k \right)^{\tt T}$
is a vector in $\C^{2k}$.
Since the random variables~\eqref{eq-20} are complex Gaussian, one can rewrite
the right-hand side of~\eqref{eq:KRH} in a more explicit form
\begin{equation}\label{eq-30}
\rho_f (z_1, ..., z_k) = \frac1{\pi^{2k} \det \Gamma_f }
\int_{\C^k} |\eta_1|^2 ... |\eta_k|^2  e^{-\frac12 \langle \Gamma_f^{-1} \eta', \eta' \rangle}
{\rm d}m(\eta_1) ... {\rm d} m(\eta_k),
\end{equation}
where $\Gamma_f = \Gamma_f (z_1, ..., z_k)$ is the covariance matrix of the random variables~\eqref{eq-20}.
We consider the linear functionals
\[
L f = \sum_{j=1}^k \left[ \alpha_j f(z_j) + \beta_j f'(z_j) \right]
= \frac1{2\pi{\rm i}} \int_{\gamma} f(z) r^L(z)\, {\rm d}z,
\]
where
\[
r^L(z) = \sum_{j=1}^k \left[ \frac{\alpha_j}{z-z_j} + \frac{\beta_j}{(z-z_j)^2} \right],
\]
and $\gamma \subset K$ is a smooth contour that bounds a domain $G'\subset K$ that contains the
points $z_1, ..., z_k$.
Observe that for every vector
$\delta = \left( \alpha_1, \beta_1, ... , \alpha_k, \beta_k \right)^{\tt T} $
in $\C^{2k}$, we have \[ \langle \Gamma_f \delta, \delta \rangle = \E | Lf |^2\,.\]
This observation allows us to estimate the matrix $\Gamma_f^{-1}$, and hence
the integral on the right-hand side of~\eqref{eq-30},
using some simple tools from the theory of analytic functions of
one complex variable.

We note that using another approach to analyzing the right-hand side of~\eqref{eq-30},
Bleher, Shiffman, and Zelditch proved in \cite{BSZ} that if the points $z_i$ are well
separated from each other, i.e.,
\[ \min_{i\ne j}|z_i-z_j|\ge \delta >0, \] then some estimate similar
to~\eqref{eq:cluster} holds with a factor $C(k, \delta)$ instead of $C_k$.
Unfortunately, in this form the result is difficult to apply. For instance, it does not
yield the boundedness of the $k$-point functions on the whole $\C^k$, and we could not
use it for the proof of Theorem~\ref{thm_normality}. On the other hand, the result of
Bleher, Shiffman and Zelditch is valid for a wider class of zero point processes.

\bigskip

\bigskip

\medskip

\centerline{\sc\Large Part~II. Random nodal lines}

\section{Gaussian spherical harmonic and Gaussian plane wave}

We introduce two remarkable Gaussian random functions closely related
to each other: the {\em Gaussian spherical harmonic} on the two-dimensional sphere
$\Ss$ and its scaling limit, the {\em Gaussian plane wave}.
The study of random plane waves, and in particular, of their nodal
portraits, originated in applied mathematics and goes back to
M.~S.~Longuet-Higgins~\cite{LH} who computed various
statistics of nodal lines for Gaussian random waves in connection
with the analysis of ocean waves. One of the reasons for the
recent interest in random plane waves
is the heuristic principle proposed by M.~V.~Berry~\cite{Berry}
called `the random wave conjecture'.
This principle says that the behaviour of high-energy Laplace
eigenfunctions in the case when the corresponding geodesic flow is ergodic
(the so called `highly excited quantum chaotic eigenfunctions')
should resemble the behaviour of Gaussian
random waves. More generally, one would expect that the random spherical
harmonic can serve as a good model for the typical behaviour of high-energy Laplace eigenfunctions
on a compact surface endowed with a smooth Riemannian metric.

\subsection{Spherical harmonics}\mbox{}

\medskip\noindent
The spherical harmonic of degree $n$ is a real-valued  eigenfunction of
the Laplacian (with the minus sign)
on the two-dimensional sphere $\Ss$ corresponding to the eigenvalue
$ \lambda_n = n(n+1) $. Equivalently, it is a trace of a homogeneous
harmonic polynomial in $\mathbb R^3$ of degree $n$ on the unit sphere.
Let $\HH_n$ be the $2n+1$-dimensional real Hilbert space of
spherical harmonics of degree $n$ equipped with the $L^2(\Ss)$-norm. The
{\em Gaussian spherical harmonic} $f$ is the sum
\[
f_n = \sum_{k=-n}^{n}\xi_k Y_k
\]
where $\xi_k$ are independent identically distributed mean zero Gaussian (real)
random variables with $\E \xi_k^2 = \frac1{2n+1}$ and $\bigl\{ Y_k \bigr\}$
is an orthonormal basis of $\HH_n$, so $ \E\|f\|^2_{L^2(\Ss)}=1$.
As a random function, $f_n$  does not
depend on the choice of the basis $\bigl\{ Y_k \bigr\}$ in $\HH_n$.
Since the scalar product in the Hilbert space $\HH_n$ is invariant under
rotations of the unit sphere,
the distribution of the random spherical harmonic $f_n$ is also
rotation invariant.
The covariance function of the Gaussian spherical harmonic equals
\[
\E \bigl\{ f_n(x)f_n(y) \bigr\}
= P_n( \cos\Theta (x, y) )
\]
where $\Theta (x, y)$ is the angle between $x$ and $y$, and $P_n $ is the Legendre
polynomial of degree $n$ normalized by $P_n(1)=1$.

\subsection{Random plane waves}\mbox{}

\medskip\par\noindent
Now, we turn to the Gaussian plane wave.
Informally speaking, it is the two-dimensional Fourier transform  of the
white noise on the unit circumference $\mathbb S^1\!\subset\!\R^2$. More formally,
we start with the Hilbert space $L^2_{\tt sym}(\mathbb S^1)$ that consists of
complex valued $L^2$-functions $\phi$ on $\mathbb S^1$ satisfying the symmetry condition
\[
\phi(-\la)=\overline{\phi(\la)}, \quad \la\in\mathbb S^1,
\]
and consider the Fourier image
of this space $\mathcal H = \mathcal F L^2_{\tt sym}(\mathbb S^1)$ with the scalar product
inherited from $ L^2_{\tt sym}(\mathbb S^1) $. The space $\HH$ consists of real-analytic
functions
\[
\Phi (x) = \int_{\mathbb S^1} e^{{\rm i}x\cdot\la} \phi(\la)\, {\rm d}m(\la)
\]
($m$ is the Lebesgue measure on $\mathbb S^1$) satisfying the Helmholtz equation
$\Delta \Phi + \Phi = 0$. The Gaussian plane wave is the sum of the
random series \[ F = \sum_k \eta_k \Phi_k \]
where $\eta_k$ are standard identically distributed independent (real) Gaussian random variables,
and $ \{ \Phi_k \} $ is an orthonormal basis in $\HH$. The series converges almost surely,
and its sum is again a real analytic function in $\R^2$ satisfying the same
Helmholtz equation. This construction does not depend on the choice of the basis
$\bigl\{ \Phi_k\bigr\}$, and the distribution of the random function $F$ is invariant with respect
to translations and rotations of the plane (since the norm in $\HH$ is translation and rotation
invariant).

Applying the Fourier transform to the standard
orthonormal basis $\bigl\{ \la^m \bigr\}_{m\in\mathbb Z}$ in $L^2(\mathbb S^1)$,
we get the functions
$ {\rm i}^m J_m(r) e^{{\rm i}m\theta} $
where $(r, \theta)$ are polar coordinates, and $J_m$ is the Bessel function of order $m$.
This yields a more explicit formula for the Gaussian plane wave:
\[
F(x) = {\rm Re\,} \sum_{m\in\mathbb Z} \zeta_m J_{|m|} (r) e^{{\rm i} m\theta},
\qquad x=(r, \theta),
\]
where $\zeta_m$ are independent identically distributed complex Gaussian random variables
with $\E |\zeta_m|^2 = 2$.

The covariance function of $F$ (which is the
same as the reproducing kernel of the space $\HH$) is given by the Bessel kernel:
\[
\E \bigl\{ F(x) F(y) \bigr\} = J_0(|x-y|)\,.
\]

\medskip
It is worth mentioning that there are other constructions of random plane `monochromatic'
waves as random linear combinations (`superpositions') of elementary plane waves
$e_\la (x) = e^{{\rm i} \la \cdot x}$. For instance, following Oravecz, Rudnick, Wigman~\cite{ORW}
and Rudnick, Wigman~\cite{RW},
one can consider `arithmetic random waves'
\[
h_N(x) = {\rm Re\,} \sum_\nu \zeta_\nu e^{2\pi {\rm i} (\nu\cdot x) }
\]
where $\zeta_\nu$ are independent identically distributed complex
Gaussian random variables with $\E |\zeta_m|^2 = 2$,
and the sum is taken over $\nu\in \mathbb Z^2$ with
$|\nu|^2=N$. This model remarkably combines analysis and probability theory
with the number theory. Its covariance function
\[
\E \bigl\{ h_N(x) h_N(y) \bigr\} = \sum_{\nu} \cos 2\pi \bigl( \nu\cdot (x-y) \bigr)
\]
has a more erratic behaviour than the covariance functions of
the Gaussian spherical harmonic and the Gaussian plane wave.

\medskip

\subsection{Random plane waves as scaling limits of random spherical
harmonics}\mbox{}

\medskip\noindent
The Gaussian plane wave $F$ is {\em a scaling limit} of the Gaussian
spherical harmonic $f_n$ when $n\to\infty$. This is a very special case
of a result of Zelditch~\cite{Zelditch1} pertaining to a wide class of
Riemannian smooth surfaces, in particular, to all real-analytic
Riemannian surfaces.

Informally, for any fixed $R$,
the restrictions of the Gaussian functions $f_n$ on spherical disks of
radius $R/n$ converge as random processes to
the restriction of $F$ on the euclidean disk of radius $R$. More formally,
we fix a point $x_0\in\Ss$, and define the random Gaussian function $F_n$
on the tangent plane $T_{x_0}\Ss$ by
\begin{equation}\label{eq:scaling}
F_n(u) = \bigl( f_n \circ \exp_{x_0} \bigr) \left( \tfrac{u}{n} \right),
\end{equation}
where $ \exp_{x_0}\colon T_{x_0}\Ss \to \Ss $ is the exponential map.
After this scaling, the covariance equals
\[
\E \bigl\{ F_n (u) F_n (v) \bigr\} = P_n \left(
\cos\Theta \left( \exp_{x_0}\left( \tfrac{u}{n} \right),
\exp_{x_0}\left( \tfrac{v}{n} \right) \right)  \right)
\]
When $n$ goes to $\infty$, the angle between the points
$\exp_{x_0}\big(\frac{u}{n}\big)$, and  $\exp_{x_0} \big(
\frac{v}{n}\big)$ on the sphere is equivalent to $|u-v|/n$
(locally uniformly in $u$ and $v$).
Then by classical Hilb's
asymptotics of the Legendre polynomials~\cite[Theorem~8.21.6]{Szego},
the scaled covariance function $ \E \big\{ F_n (u)
F_n (v) \big\} $ converges to the Bessel kernel $J_0 (|u-v|)$
locally uniformly in $u$ and $v$.

\section{Nodal portrait}

In most cases, the basic questions about the asymptotic behaviour of the nodal portrait of
the Gaussian spherical harmonic $f_n$ as $n\to\infty$, and their counterparts
for the Gaussian plane wave in the `large area limit' are equivalent
to each other. In what follows, we concentrate on spherical harmonic versions
which are somewhat easier to formulate.

\medskip
For the spherical harmonic $g\in\HH_n$, we denote by
$Z(g)=\{x\in\Ss\colon g(x)=0\}$ its nodal set.
The connected components of the complement
$\Ss \setminus Z(g) $ are called nodal domains of $g$. The following (deterministic)
facts are special cases of well-known results valid for Laplace eigenfunctions
on smooth Riemannian surfaces:
\begin{theorem}
There is a positive numerical constants $C$ such that for each $g\in\HH_n$,
the nodal set $Z(g)$ is a $Cn^{-1}$-net on $\Ss$.
\end{theorem}
\begin{theorem}\label{area-estimate}
There is a positive numerical constant $c>0$ such
that for each $g\in\HH_n$, every nodal domain of $g$ contains a disk of radius
$cn^{-1}$.
\end{theorem}
Together with Figure~3, this gives a very rough idea of how
the nodal portraits of a spherical harmonic of large degree should look.
\begin{figure}[h]\label{fig:3}
\centering
\scalebox{0.45}{\includegraphics{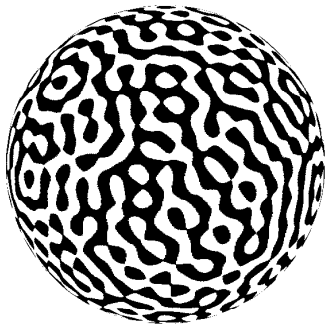}}
\caption{Nodal portrait of the Gaussian spherical harmonic of degree
$40$ (figure by A.~Barnett)}
\end{figure}

One can find more information about the geometry and the topology of the nodal portraits of
spherical harmonics (and more generally, of high-energy
Laplace eigenfunctions on smooth Riemannian surfaces) in
the pioneering works of Donnelly and Fefferman~\cite{DF1, DF2, DF3},
as well as in the more recent works of Eremenko, Jackobson, and Nadirashvili~\cite{EJN},
Mangoubi~\cite{Mangoubi},
and Nazarov, Polterovich and Sodin~\cite{NPS}. Still, our understanding of nodal portraits
is rather restricted, and, in our opinion, this classical area of analysis is very much underdeveloped.

\subsection{Length of the nodal set}\mbox{}\nopagebreak

\medskip\noindent
The basic characteristics of the nodal set of a spherical harmonic $g$
are its length $ L(g) $ and the number  $ N(g) $ of connected components
(which is one less than the number of nodal domains).
Useful classical integral formulas for the length due to Poincar\'e
and to Kac and Rice make the length a somewhat easier object for a study.
For instance, one can prove
\begin{theorem}\label{length-estimate}
There exists a positive numerical constant $C$ such that for each $g\in\HH_n$,
$ C^{-1} n \le L(g) \le Cn $
\end{theorem}
This is a special case of a more general result valid for Laplace
eigenfunctions corresponding to large eigenvalues (with $n$ replaced by $\sqrt{\la}$).
The lower bound is valid for any smooth Riemannian surface (this is a result
of Br\"uning~\cite{Br}), while the upper bound was proven by Donnelly and Fefferman~\cite{DF1}
for real-analytic surfaces.
In the smooth category, it was conjectured by S.~T.~Yau,
and still remains open in spite of many efforts.
Note that one can easily deduce the upper bound in Theorem~\ref{length-estimate}
from the fact that spherical harmonics are restrictions of polynomials (that is,
without using the deep result of Donnelly and Fefferman).

\medskip
For the Gaussian spherical harmonic, B\'erard showed in~\cite{Berard} that
\begin{theorem} $\E L(f_n) = \pi  \sqrt{2\la_n} = \sqrt{2}\,\pi n + O(1)$.
\end{theorem}

The question about the variance is more delicate.
Recently, Wigman~\cite{Wigman} confirmed a guess made
by M.~V.~Berry~\cite{Berry2} in a slightly different context:
\begin{theorem}
For $n\to\infty$,
 \[
 {\rm variance\ of\ } L(f_n) = \frac{65}{32} \log n + O(1)\,.
 \]
\end{theorem}
\par\noindent
The proof of this theorem is based on a very careful analysis of asymptotic cancelations
that appear in the Kac-Rice integral representation  of the variance of $L(f_n)$.

\subsection{The number of connected components}\mbox{}

\medskip\par\noindent
There are few classical facts about the number of components $ N(g) $.
The celebrated Courant nodal domain theorem yields
\begin{theorem} For every $g\in\HH_n$, $N(g)\le n^2$.
\end{theorem}
\par\noindent For large $n$, this upper bound was improved by Pleijel~\cite{Pl} to
$ 0.69 n^2$. Apparently, the sharp asymptotic upper bound is not known yet. Simple examples
show that it cannot be less than $(\tfrac 12 + o(1))n^2$.
H.~Lewy~\cite{Lewy} gave an elegant
construction of spherical harmonics of any degree $ n $ whose nodal sets have
one component for odd $ n  $ and two components for even $ n $, which proves
that no non-trivial lower bound for $N(g)$ is possible.

Till recently, nothing had been known about the asymptotic
properties of the random variable $N(f_n)$ when the degree $n$ is large.
The principal difficulty is its {\em non-locality}:
observing the nodal curves only locally, one cannot make any definite conclusion about the
number of connected components.
Several years ago Blum, Gnutzmann, and Smilansky~\cite{BGS} raised
a question about the distribution of the number of nodal domains
of high-energy Laplace eigenfunctions. In the ergodic case, in accordance
with Berry's heuristic principle, they suggested to find
this distribution for Gaussian random plane waves and performed
the corresponding numerics. To compute this distribution,
Bogomolny and Schmit proposed in~\cite{BS} an elegant
percolation-like lattice model for description of nodal domains of
random Gaussian plane waves. This model completely ignores the (quite big)
correlations between the values of the random function $f_n$ at different
points but nevertheless agrees with numerics pretty well.
This agreement is probably due to some hidden `universality law' rather
then the possibility to directly reduce one model to another.

\subsection{Bogomolny-Schmit percolation-like model}\mbox{}

\medskip\noindent
The Bogomolny-Schmit {\em hypothesis} is that the distribution of nodal domains $ N(f_n)$ is
roughly the same as in the following critical percolation model.
Consider the square lattice with the total number
of sites equal to $\bigl( \E L(f_n) \bigr)^2$, that is proportional to $n^2$, and change at each site
the line crossing to one of the two equiprobable avoided crossing, as shown in the following figure.
\begin{figure}[h]
\centering
\includegraphics{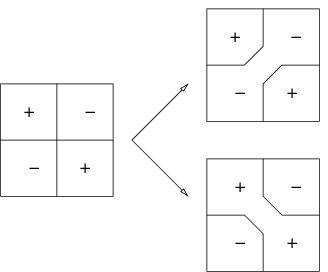}
\caption{Avoided nodal crossings in the Bogomolny-Schmit model}
\end{figure}
At different sites, the changes are independent.

Then Bogomolny and Schmit introduce two dual square lattices: the `blue one' with
vertices at the cells of the grid where the function is positive, and the `red one'
with vertices at the cells of the grid where the function is negative. Each realization
of the random process generates two graphs, the blue one whose vertices are the blue
lattice points and the red one whose vertices are the red lattice points.
Two vertices are connected by an edge if the corresponding cells of the grid
belong to the same nodal domain of the random function.
\begin{figure}[h]
\centering
\includegraphics{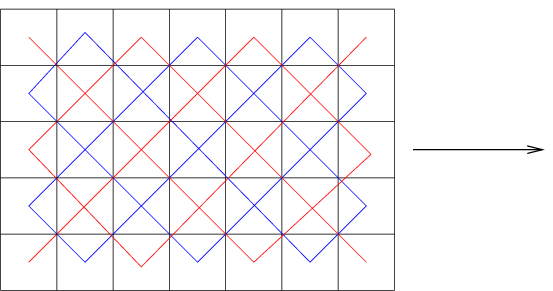}\hspace{0.12in}
\includegraphics{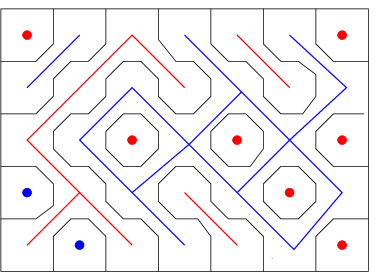}
\caption{Bond percolation on two dual lattices}
\end{figure}
Each of these graphs uniquely
determines the whole picture, so it suffices to consider only one of them,
and each of them represents the critical bond percolation on the
corresponding square lattice.
Then using some heuristics
coming from statistical mechanics, Bogomolny and Schmit predicted that for $n\to\infty$,
\[
\E N(f_n) = ( a+o(1) ) n^2,
\]
and
\[
{\rm variance\ of\ } N(f_n) = ( b+ o(1) )n^2,
\]
with explicitly computed  positive numerical constants $a$ and $b$. They also argued that the
fluctuations of the random variable $ N(f_n) $ are asymptotically Gaussian when $n\to\infty$,
and concluded their work with a remarkable prediction of the power distribution law for
the areas of nodal domains, based on the percolation
theory.

It would be interesting to test numerically
whether the Bogomolny-Schmit model persists
for random linear combinations of plane waves $e^{{\rm i}k\cdot x}$
with {\em different} wave numbers $k$.

\subsection{Rigorous results}\mbox{}

\medskip\noindent
Recently, we showed in~\cite{Nazarov-Sodin} that, in accordance with one of
the Bogomolny and Schmit predictions, $\E N(f)/n^2$ tends to a positive
limit when $n\to\infty$,
though our proof does not provide us with an explicit value of the limit
$a$, so we cannot juxtapose it with the one predicted by Bogomolny and Schmit.
In addition, we proved that the random variable
$N(f)/n^2$ concentrates around this limit exponentially.
Since for any spherical harmonic $g\in\HH$,
the total length of its nodal set $Z(g)$ does not exceed ${\rm Const}\, n$,
our result yields that, for a typical spherical
harmonic, most of its nodal domains have diameters comparable to
$1/n$.

\begin{theorem}[Number of nodal domains]\label{thm.main}
There exists a constant $a>0$ such that, for every $\epsilon > 0$, we have
$$
\PP \left\{\left|\frac{N(f_n)}{n^2}-a\right|>\e\right\}\le
C(\e)e^{-c(\e)n}
$$
where $c(\e)$ and $C(\e)$ are some positive constants depending on
$\e$ only.
\end{theorem}

The exponential decay in $n$ in Theorem~\ref{thm.main} cannot
be improved: we showed that given a positive and arbitrarily small $\kappa$,
\[
\PP\left\{ N(f_n) < \kappa n^2 \right\} \ge e^{-C(\kappa)n}\,.
\]
On the other hand, our proof of Theorem~\ref{thm.main} gives a
very small value $c(\epsilon) \simeq \epsilon^{15}$ and it would
be nice to reduce the power $15$ of $\epsilon$ to something more
reasonable. The question about the variance of $N(f_n)$ remains open.

\medskip The last but not least remark is
that the proof of Theorem~\ref{thm.main} uses only
relatively simple tools from the classical analysis, which we
believe may work in a more general setting of random functions of
several real variables (and for higher Betti numbers), while
it seems that the Bogomolny-Schmit
model is essentially a two-dimensional one.

\subsection{Related work}\mbox{}

\medskip
We are aware of several encouraging attempts to tackle
similar questions in different contexts.
In~\cite{Swerling} (motivated by some engineering problems),
Swerling estimated from below and from above the mean number of
connected components of the {\em level lines} $Z(t, f) = \{ f=t\}$
of a random Gaussian trigonometric polynomial $f$ of two variables
of a given degree $n$. His method is based on estimates of the
integral curvature of the level line $Z(t, f)$. The estimates are
rather good when the level $t$ is separated from zero, but as
$t\to 0$ they are getting worse and, unfortunately, give nothing
when $t=0$.

\medskip
In the paper~\cite{Malevich}, Malevich considered $C^2$-smooth Gaussian
random functions $f$ on $\mathbb R^2$ with positive covariance function
that decays polynomially as the distance between the points tends to infinity.
She proved
that for $T\ge T_0$,
\[
C^{-1} T^2 \le \E N(T) \le C T^2\,,
\]
where $N(T)$ is the number of the connected components of the zero set
of $f$ that are contained in the square $[0, T] \times [0, T]$,
and $C$ is a positive numerical constant. Her proof relies heavily on
the positivity property of the covariance function that does not hold for
Gaussian spherical harmonics or for Gaussian trigonometric polynomials.

\medskip
In the recent paper~\cite{MW}, Mischaikow and Wanner
studied the following question. Suppose $f$ is a random smooth function on
the square $[0, 1]^2$ with periodic boundary conditions and that the signs
of $f$ are computed at the vertices of the grid with mesh $\delta$.
{\em How small must $\delta$ be (in terms of the a priori smoothness constants of $f$) in order
to recover the Betti numbers of the sets $\{f>0\}$
and $\{f<0\}$ with probability close to one ?}
In particular, they show that for random trigonometric polynomials
of two variables of degree $N$, it suffices to take $\delta=cN^{-2}$ where $c$ is a
sufficiently small positive numerical constant.
It is possible that their bounds can be significantly improved
if instead of recovering the exact values of the Betti numbers one tries
to recover them with a small relative error.

\section{The sketch of the proof of the theorem on the number of nodal domains}

Here, we will describe the main ideas behind the proof of Theorem~\ref{thm.main}. All the details
can be found in~\cite{Nazarov-Sodin}.

\subsection{The lower bound $\E N(f_n) \ge cn^2$}\mbox{}

\medskip\noindent
This is the simplest part of the story. Denote by $d(.,.)$ the spherical distance.
Given a point $x\in\Ss$ and a large positive constant $C$, we consider the event
\[
\Omega_x = \bigl\{ f_n(x) > C,
\ {\rm and }\ f_n(y)<-C {\rm\ for\ all}\ y {\rm\ satisfying}\  d(x,y)=\frac{\rho}n \bigr\}\,,
\]
where $\rho$ is a constant whose value will be specified below.
Clearly, if the event $\Omega_x$ occurs, then the disk of radius $\rho/n$ centered at $x$
contains a closed nodal line of $f_n$.
We claim that
\[
\PP (\Omega_x) \ge c > 0\,,
\]
where $c$ is a positive constant. The reason is rather straightforward:
for every point $x\in\Ss$,
there exists a function $b_x\in\HH_n$ with $\| b_x\|=1$ such that
\[
b_x(x) > c_0 \sqrt{n} \ \ {\rm and }\ b_x (y) < -c_0\sqrt{n} {\rm\ whenever\ }
d(x,y)=\frac{\rho}{n}\,.
\]
One can take as
$b_x$ the zonal spherical harmonic with ``pole'' $x$.
Then we can represent $f_n$ in the form
\[
f_n = \xi_0 b_x + f_x
\]
where $\xi_0$ is a Gaussian random variable with $\E \xi_0^2 = \frac1{2n+1}$, and $f_x$ is
a Gaussian spherical harmonic with $\E \|f_x\|^2 = \frac{2n}{2n+1} $ independent of $\xi_0$,
and check that with positive probability, the `perturbation' $f_x$ cannot destroy a short nodal curve around
point $x$ provided by the function $b_x$.

It remains to pack the sphere $\Ss$ by $ \simeq n^2 $ disjoint disks of radius $2\rho/n$. With a positive
probability, each of these disks contains a closed nodal line of $f_n$. Whence, the lower bound for $\E N(f_n)$.

\subsection{Levy's concentration of measure principle}\mbox{}

\medskip\noindent
To establish the exponential concentration of the random variable $N(f_n)$ around its
median, we would like to use a version of classical Levy's concentration of measure
principle.

Given a set $K$, we denote by $K_{+\rho}$ the $\rho$-neighbourhood of
$K$. We apply this notation to subsets of $\HH_n$ and the
$L^2$-distance, to subsets of $\Ss$ and the usual spherical
distance, and also to subsets of $\R^d$ with the Euclidean distance.
The following Gaussian isoperimetric theorem is due to Sudakov and Tsirelson~\cite{ST}
and Borell~\cite{Borell}:
\begin{theorem}\label{thm:isoperim}
Let $\gamma_d$ be the standard Gaussian measure on $\R^d$.
Let $\Sigma\subset\R^d$ be a Borel set, and $\Pi$ be an affine half-space such that
\[
\gamma_d (\Sigma) = \gamma_d (\Pi)\,.
\]
Then for each $t>0$,
\[
\gamma_d (\Sigma_{+\rho}) \ge \gamma_d (\Pi_{+\rho})\,.
\]
\end{theorem}
A simple computation shows that if $\gamma_d (\Pi_{+\rho})$ is not too
close to $1$, then $\gamma_d (\Pi)$ must be exponentially small
in $d$, like $\exp [-c\rho^2 d]$. Applying this to the $2n+1$-dimensional
space $\HH_n$ of spherical harmonics of degree $n$, we get
\begin{corollary}[Concentration of Gaussian measure on $\HH_n$]\label{lemma:concentr}
Let $G\!\subset\!\HH_n$ be any measurable set of spherical harmonics.
Suppose that the set $G_{+\rho}$ satisfies $\PP(G_{+\rho})<\frac 34$.
Then $\PP(G)\le 2 e^{-c\rho^2 n}$.
\end{corollary}

To use the concentration of measure principle, we need to show that the
number $N(f)$ doesn't change too much under slight perturbations of $f$
in the $L^2(\Ss)$-norm. Certainly, this is not true for all $f\in \HH_n$,
but we will show that the
``unstable'' spherical harmonics $f\in\HH_n$ for which small perturbations
can lead to a drastic decrease in the number of nodal lines are
exponentially rare. Here is a key lemma which is probably the most novel
part of the whole story:
\begin{lemma}[Uniform lower semi-continuity of $N(f_n)/n^2$]\label{lemma}
For every $\epsilon > 0$, there exist $\rho>0$ and an exceptional set
$E\subset \HH_n$ of probability $\PP(E)\le C(\e)e^{-c(\e)n}$ such
that for all $f\in \HH_n\setminus E$ and for all $g\in \HH_n$
satisfying $\|g\| \le \rho$, we have \[ N(f+g)\ge N(f)-\e n^2\,. \]
\end{lemma}

\smallskip
The uniform lower semi-continuity lemma readily yields the exponential concentration of
the random variable $N(f_n)/n^2$ near its median $a_n$. First, consider the set
\[ G = \bigl\{ f\in\HH_n \colon N(f)>(a_n+\epsilon)n^2 \bigr\}.\]
Then for $ f\in (G\setminus E)_{+\rho}$, we have $N(f)>a_nn^2$,
and therefore, $ \PP \bigl( (G\setminus E)_{+\rho} \bigr) \le \frac 12 $. Hence,
by the concentration of Gaussian measure, $\PP(G\setminus E)\le 2e^{-c\rho^2 n}$, and
finally,
\[
\PP(G)\le \PP(G\setminus E) + \PP (E) \le 2 e^{-c\rho^2n}+C(\e)e^{-c(\e)n}\le C(\e)e^{-c(\e)n}\,.
\]
Now, we turn to the set
\[ F = \bigl\{ f\in \HH_n \colon N_f<(a_n-\epsilon)n^2 \bigr\}.\]
Then
\[ F_{+\rho} \subset \bigl\{ f\in\HH_n\colon  N_f < a_n n^2\bigr\} \cup E\,,\]
so that
\[
\PP (F_{+\rho}) \le \frac12 + C(\e)e^{-c(\e)n} < \frac34
\]
for large $n$, and it follows that $ \PP(F)\le 2e^{-c\rho^2 n} $.

\subsection{The uniform lower continuity of the functional $f\!\mapsto\!N(f)$ outside of
an exceptional set}\mbox{}

\medskip\par\noindent Here we explain how we prove Lemma~\ref{lemma}.

\subsubsection{Exceptional spherical harmonics $E$ with unstable nodal portraits}
Instability of the nodal portrait of a spherical harmonic $f\in\HH_n$ under small perturbations
is caused by points where $f$ and $\nabla f$ are simultaneously small. Let $\alpha$ and $\delta$ be
small positive parameters, and let $R$ be a large positive parameter (all of them will depend on
$\epsilon$ from Lemma~\ref{lemma}).
Cover the sphere $\Ss$ by $ \simeq R^{-2} n^2 $ disks $D_j$ of radius $R/n$
in such a way that the concentric disks $4D_j$ with $4$ times larger radius
cover the sphere with a bounded multiplicity.
We call the disk $D_j$ {\em stable} if for each $x\in 3D_j$ either
$ |f(x)|\ge \alpha $ or $ |\nabla f(x)| \ge \alpha n $. Otherwise, the disk $D_j$ is {\em unstable}.
We call the spherical harmonic $f\in\HH_n$ {\em exceptional} if the number of unstable disks
is at least $\delta n^2$, and denote by $E$ the set of all exceptional spherical harmonics
of degree $n$.
\begin{lemma}\label{lemma:exc-set}
Given $\delta >0$,
there exist positive $C(\delta)$ and $c(\delta)$ such that
\[
\PP (E) \le C(\delta) e^{-c(\delta)n}
\]
provided that the constant $\alpha$ is sufficiently small.
\end{lemma}
\par\noindent
Curiously, the proof of this lemma uses the concentration of measure
principle again. It also uses the fact that given $x\in\Ss$, the Gaussian random
variable $f(x)$ and the Gaussian random vector $\nabla f(x)$ are independent.

\subsubsection{Identification of unstable connected components}
It remains to show that at most $\epsilon n^2$ nodal components of
a stable spherical harmonic can disappear after perturbation of $f$ by another
spherical harmonic $g\in\HH_n$ with sufficiently small $L^2$-norm.
First, in several steps, we identify possibly `unstable' connected
components of the nodal set $Z(f)$ that can disappear after perturbation,
show that their number is small compared to $n^2$, and discard them.
Then we verify that all other connected components of $Z(f)$ do not disappear
after the perturbation.

\smallskip\noindent\underline{First},
we discard the nodal components $\Gamma$ whose diameter is bigger than $R/n$.
By the upper bound in the length estimate in Theorem~\ref{length-estimate},
their number is $ \lesssim CR^{-1}n^2 $
which is small compared to $n^2$.

With each remaining component $\Gamma$ of the nodal
set $Z(f)$ we associate a disk $D_j$ such that $D_j \cap \Gamma \neq \emptyset$.
Then $\Gamma\subset 2D_j$. Since each nodal domain contains a disk of radius
$c/n$ (Theorem~\ref{area-estimate}), the number of components $\Gamma$ intersecting
$D_j$ (and, thereby, contained in $2D_j$) is bounded.

\smallskip\par\noindent\underline{Second}, we discard the components $\Gamma$ with unstable
disks $D_j$. Since $f$ is not exceptional, and since
each disk $D_j$ cannot intersect too many components contained in $2D_j$,
the number of such components is also small compared to $n^2$.

\smallskip\par\noindent\underline{At last},
we discard the components $\Gamma$ such that
\[
\max_{3D_j} |g| \ge \alpha\,.
\]
To estimate the number $N$ of such disks, we denote by $D_j^* \subset 4D_j$
the disk of radius $1/n$ centered
at the point $y_j$ where $|g|$ attains its maximum in $3D_j$.
By standard elliptic estimates,
\[
\int_{D_j^*} |g|^2 \gtrsim n^{-2} |g(y_j)| = \alpha^2 n^{-2}\,,
\]
whence
\[
\rho^2 \ge \| g \|_{L^2(\Ss)} \gtrsim  N \alpha^2 n^{-2}\,,
\]
that is, $ N \lesssim \rho^2 \alpha^{-2} n^2 $. As above, we conclude that
the number of components $\Gamma$ affected by this is
$ \lesssim  R^2N \lesssim R^2 \rho^2 \alpha^{-2} n^2 $
which is much less than $\epsilon n^2$ provided that $\rho^2 $ is much less
than $ \epsilon \alpha^2 R^{-2}$.

\subsubsection{Verification of stability of the remaining connected components}
Now, we claim that the remaining components $\Gamma$ cannot be affected by
the perturbation of $f$ by $g$. To see this, we consider the connected component
$A_\Gamma (t)$ of the set $\{ |f|<t \} $ that contains
$\Gamma$, and look what may happen with this component when $t$ grows from $0$ to
$\alpha$.
\begin{itemize}
\item[$\bullet$]
As long as $A_\Gamma (t)$ stays away from the
boundary $\partial (3D_j)$, it cannot merge with another
component of $\{|f|<t\}$  because such a merge can occur only at a
critical point of $f$ and there are none of them in $A_\Gamma (t) \cap 3D_j$.
\item[$\bullet$]
For the same reason,
neither of the two boundary curves of $A_\Gamma (t)$ can collapse
and disappear.
\item[$\bullet$]
At last, $A_\Gamma (t)$ cannot reach
$\partial (3D_j)$ before it merges with some other component:
indeed, if $x\in A_\Gamma (t)$ and $A_\Gamma (t)$ lies
at a positive distance from the boundary $\partial (3D_j)$
then we can go from $x$ in the direction of $\nabla f$ if $f(x)<0$
and in the direction of $-\nabla f$ if $f(x)>0$. In any case, since
$|\nabla f|>\alpha n$ in $A_\Gamma (t)$, we shall reach the zero set
$Z(f)$ after going the length  $1/n$ or less. Since $\Gamma$ is the only
component of $Z(f)$ in $A_\Gamma (t)$ before any merges,
we conclude that $A_\Gamma (t) \subset \Gamma_{+1/n}$.
Recalling that ${\rm dist}\, (\Gamma, \partial (3D_j)\,)>R/n$, we
see that, for each $t\le \alpha$, $A_\Gamma (t)$ stays away from the
boundary $\partial (3D_j)$.
\end{itemize}

Thus, each component $\Gamma$ lies in a topological annulus $A_\Gamma =
A_\Gamma (\alpha)$ which is contained with its boundary in the open
disk $3D_j$ and such that $f=+\alpha$ in one boundary curve of
$A_\Gamma$ and $f=-\alpha$ on the other. Recalling that $|g|<\alpha$ in
$3D_j$, we conclude that $Z(f+g)$ has at least one
connected component in $A_\Gamma$.

\subsection{Existence of the limit of $ \E N(f_n) / n^2$ }\mbox{}

\medskip\noindent
We already know that $\E N(f_n) \gtrsim n^2$ and that
$N(f_n)/n^2$ concentrates near its median exponentially. Thus,
to finish the proof of Theorem~\ref{thm.main}, it remains to
show that the sequence $ \bigl\{ \E N(f_n)/n^2 \bigr\} $ converges.
We deduce this from the fact that the Gaussian spherical harmonic
has a scaling limit combined with rotation invariance of the distribution
of $f_n$.
Since this part does not require any new ideas beyond
the ones we've already introduced, we just refer the reader to~\cite{Nazarov-Sodin}
for the details.

\subsection{Comments and questions}\mbox{}

\medskip\par\noindent
Making use of a non-critical version of their percolation model, Bogomolny and Schmit
obtained in~\cite{BS2} a series of predictions for the behaviour of the
components of level sets which agree with numerics.
In a stark contrast, we do not have a rigorous answer even to
the following most basic question:
\begin{question}
Prove that for each $\epsilon\!>\!0$ and each $\eta>0$, the probability
that the level set $\bigl\{x\in\Ss\colon  f_n(x)\!>\! \epsilon \bigr\}$ has a component
of diameter larger than $\eta$ tends to zero as $n\to\infty$.
\end{question}

One of the reasons for our ignorance is the aforementioned non-locality of the
number of connected components.
Another essential difficulty is a very slow decay of the correlations
which does not allow us to think of restrictions of
our process to a collection of well-separated disks as of almost independent processes.

\begin{question}
Estimate the mean number
of large components of the nodal set whose diameter
is much bigger than $1/n$. For instance, of those whose diameter
is comparable to $n^{-\alpha}$ with $0<\alpha<1$.
\end{question}

\medskip
Nothing is known about the number of connected components of the nodal set
for `randomly chosen' high-energy Laplace eigenfunction $f_\la$  on an
arbitrary compact surface $M$ without boundary endowed with a smooth
Riemannian metric $g$.
It is tempting to expect that Theorem~\ref{thm.main} models what is happening
when $M$ is the two-dimensional sphere $\Ss$ endowed with a generic Riemannian metric
$g$ that is sufficiently close (with several derivatives) to the constant one.

Instead of perturbing the `round metric' on the sphere $\Ss$,
one can add a small potential $V$ to the Laplacian on the sphere.
The question remains just as hard.

\subsection*{Acknowledgements}\mbox{}

\medskip\par\noindent
We are grateful to Manjunath Krishnapur, Yuri Makarychev, Yuval Peres,
Leonid Polterovich, Ze\'ev Rudnick, Bernard Shiffman, Boris Tsirelson, Sasha Volberg,
and Steve Zelditch for many helpful conversations
on the subject of these notes, and
to Alex Barnett, Manjunath Krishnapur, and Balint Vir\'ag  for providing us
with inspiring computer generated pictures.

\end{document}